\documentclass[11pt]{amsart}

\usepackage{amssymb}
\usepackage{amsfonts,amsmath,amssymb,mathrsfs,verbatim}
\usepackage[dvips]{graphicx}
\usepackage{graphics}
\usepackage{psfrag}
\usepackage{hyperref}
\usepackage{marginnote}
\usepackage{xcolor}
\usepackage[dvips]{graphicx}
\usepackage{pgfplots}
\pgfplotsset{compat=1.16}
\usetikzlibrary{arrows}
\usepackage{tikz}
  \usetikzlibrary{shapes,backgrounds}
  \usetikzlibrary{calc}
  \usetikzlibrary{matrix}
\usepackage{caption}
\usepackage[margin=3cm]{geometry}
\usepackage{multicol}

\usepackage{todonotes}

\newcommand{\N}{\mathbb{N}}

\newcommand{\Z}{\mathbb{Z}}

\newcommand{\lcm}{\mathrm{lcm}}

\newcommand{\mS}{\mathcal{S}}

\newcommand{\bx}{\textbf{x}}
\newcommand{\lS}{\leq_{\mS}}
\newcommand{\x}{\mathbf{x}}

\newcommand{\betti}[1]{\operatorname{Betti}(#1)}
\newcommand{\bmin}[1]{\operatorname{Minimals}_{\le \mS}(\operatorname{Betti}(#1))}

\newtheorem{thm1}{Theorem}[section]
\newtheorem{lem1}[thm1]{Lemma}
\newtheorem{cor1}[thm1]{Corollary}

\newtheorem{prop1}[thm1]{Proposition}
\newtheorem{conj1}[thm1]{Conjecture}
\newtheorem{quest1}[thm1]{Question}

\theoremstyle{definition}
\newtheorem{ex1}[thm1]{Example}
\newtheorem{def1}[thm1]{Definition}

\theoremstyle{remark}
\newtheorem{rem1}[thm1]{Remark}

\subjclass{20M14, 20M05, 13F65, 13P10}
\keywords{free numerical semigroup, universally free numerical semigroup, Betti divisible, Markov bases, universal Markov basis, circuits, universal Gr\"obner basis, Graver basis, toric ideals}

\begin{document}

\title{Universally free numerical semigroups}

\author[I. Garc\'ia-Marco]{Ignacio Garc\'ia-Marco}
\address{Facultad de Ciencias and Instituto de Matem\'aticas y Aplicaciones (IMAULL), Secci\'on de Matem\'aticas, Universidad de La Laguna, Apartado de Correos 456, 38200 La Laguna, Spain}
\email{iggarcia@ull.edu.es}
\author[P. A. Garc\'ia-S\'anchez]{Pedro A. Garc\'ia-S\'anchez}
\address{Facultad de Ciencias, Universidad de Granada, 18010 Granada, Spain}
\email{pedro@ugr.es}
\author[I. Ojeda]{Ignacio Ojeda}
\address{Facultad de Ciencias, Universidad de Extremadura, 06071 Badajoz, Spain}
\email{ojedamc@unex.es}
\author[Ch. Tatakis]{Christos Tatakis}
\address{Department of Mathematics, Aristotle University of Thessaloniki, 54124 Thessaloniki, Greece}
\email{chtatakis@math.auth.gr}

\begin{abstract}
A numerical semigroup is said to be universally free if it is free for any possible arrangement of its minimal generating set. In this work, we establish that toric ideals associated with universally free numerical semigroups can be generated by their set of circuits. Additionally, we provide a characterization of universally free numerical semigroups in terms of Gr\"obner bases. Specifically, a numerical semigroup is universally free if and only if all initial ideals of its corresponding toric ideal are complete intersections. Furthermore, we establish several equalities among the toric bases of a universally free numerical semigroup.

We provide a complete characterization of $3$-generated universally free numerical semigroups in terms of their minimal generating sets, and by proving the equality of certain toric bases. We compute exactly all the toric bases of a toric ideal defined by a 3-generated universally free numerical semigroup. Notably, we answer some questions posed by Tatakis and Thoma by demonstrating that toric ideals defined by $3$-generated universally free numerical semigroups have a set of circuits and a universal Gr\"obner basis of size 3, while the universal Markov basis and the Graver basis can be arbitrarily large. We present partial results and propose several conjectures regarding universally free numerical semigroups with more than three generators.

\end{abstract}

\maketitle

\section{Introduction}

Set $A:=\{a_1,\ldots,a_n\} \subseteq \mathbb{Z}^+$ and let \[ \mS := \langle A \rangle = \langle a_1,\ldots,a_n \rangle = \{u_1 a_1+\cdots+ u_n a_n \ | \ u_i \in\mathbb{N}\}\] be the submonoid of $\mathbb N$ generated by $A$.  Let $\mathbb K[\x] := \mathbb{K}[x_1,\ldots,x_n]$ be the polynomial ring in $n$ indeterminates over an arbitrary field $\mathbb{K}$. For $\mathbf u=(u_1,\ldots,u_n) \in \mathbb{N}^n$, we define the $A$-\emph{degree} of the monomial $\bx^{\mathbf u}:=x_1^{u_1} \cdots x_n^{u_n}$ to be $\deg_{A}(\bx^{\mathbf u}):=u_1 a_1+\cdots+ u_n a_n \in \mS$. The toric ideal of $A$, denoted $I_A$, is the binomial prime ideal of height $n-1$ generated by all the pure difference $A$-homogeneous binomials in $\mathbb K[\x]$, that is, \[ I_A = \langle \bx^{\mathbf u}- \bx^{\mathbf v} \mid \deg_{A}(\bx^{\mathbf u})=\deg_{A}(\bx^{\mathbf v}) \rangle. \]  We refer the reader to \cite{CLS, ES, ST,V} for a detailed study of toric ideals. 

Note that if $A' := A/\gcd(A) = \{a_1/\gcd(A),\ldots,a_n/\gcd(A)\}$, then $I_A = I_{A'}$.  Therefore, it can be assumed without loss of generality that the set $A = \{a_1,\ldots,a_n\}$ consists of relatively prime positive integers. In this case, $\mS$ is said to be a \emph{numerical semigroup}. For more information on submonoids of $\N$ and numerical semigroups, see \cite{AG,RA}.

Since $I_A$ has height $n-1$, we have that $I_A$ is a complete intersection (or, equivalently, $\mS$ is a complete intersection) if and only if it can be generated by $n-1$ polynomials, which, by the way, can be taken to be pure difference $A-$homogeneous binomials. Complete intersection numerical semigroups have been extensively studied in the literature, see, for instance, \cite{AG2,BGS,D,DMS,H,Stam}.

For a monomial order $\prec$ on $\mathbb{K}[\bx]$, we denote by  $\operatorname{in}_{\prec}(I_A)$ the initial ideal of $I_A$ with respect to $\prec$. Since $\operatorname{ht}(\operatorname{in}_{\prec}(I_A)) = \operatorname{ht}(I_A) = n-1$, the ideal $\operatorname{in}_{\prec}(I_A)$ is a complete intersection if and only if it can be generated by $n-1$ monomials. In other words, $\operatorname{in}_{\prec}(I_A)$ is a complete intersection if and only if the reduced Gr\"obner basis of $I_A$ with respect to $\prec$ consists of $n-1$ binomials. Furthermore, since the Gr\"obner bases are generating sets of the ideal, if $\operatorname{in}_{\prec}(I_A)$ is a complete intersection for some monomial order $\prec$ on $\mathbb{K}[\bx]$, then $I_A$ is too. By \cite[Theorem 4.7]{GT}, $I_A$ has a complete intersection initial ideal if and only if $\mS$ is free for an arrangement of its minimal generators (see Definition \ref{freesemigroups}), which makes free numerical semigroups a remarkable subfamily of complete intersection numerical semigroups.


This paper introduces and studies from different perspectives the notion of a universally free numerical semigroup. A numerical semigroup is said to be universally free if it is free for any arrangement of its minimal generating set. Although the leitmotiv of this manuscript is the concept of universally free numerical semigroup, we study other families of numerical semigroups and how the toric bases of their corresponding defining ideals are related. We have a clear picture on how these families are related for embedding dimension three, but for higher embedding dimension things are much more complicated.

The manuscript is organized as follows.

In Section \ref{Section 2} we present some fundamental facts related to several toric bases of toric ideals. The toric bases considered include the Graver basis, the universal Markov basis, the universal Gr\"obner basis, the set of circuits and the set of critical binomials of a toric ideal. 

Section 3 is divided into four subsections and is devoted to the study of different families of numerical semigroups. In the first subsection, we recall several results concerning free semigroups. The second subsection is focused on universally free numerical semigroups. We establish in Theorem \ref{universally-Grobner} that universally free numerical semigroups satisfy (\ref{eq:muin}). Specifically, we prove that for $A \subset \Z^+$, all initial ideals of $I_A$ are complete intersections if and only if $\mS$ is universally free. Additionally, we demonstrate in Proposition \ref{prop:free} that whenever $\mS$ is universally free, its toric ideal can be generated by circuits. 
The third subsection examines Betti divisible numerical semigroups, a notable subfamily of the universally free ones. In this context we establish some equalities between toric bases. The last subsection concerns circuit numerical semigroups. In Proposition \ref{pr:circuit} we provide a family of circuit numerical semigroups and, in Proposition \ref{pr:em3circuit} we prove that for $n = 3$ all circuit numerical semigroups are of the form described in Proposition \ref{pr:circuit}. Propositions \ref{prop:free}, \ref{pr:circuit} and \ref{pr:em3circuit} contribute to \cite[Open problem 5.4]{GT} where it is asked to characterize all circuit numerical semigroups.

In Section 4, we examine the case where $n = 3$, providing complete descriptions of the universally free numerical semigroups both in terms of the values of $a_1,a_2,a_3$ and computing exactly all the toric bases of the ideal $I_A$. As a consequence, by applying the above characterizations, we are able to answer an open problem posed in \cite{TTbases}, in which the authors asked to prove that the size of the Graver basis of a toric ideal cannot be bounded above by a polynomial expression on the size of the universal Gr\"obner basis or the set of the circuits of the ideal, see Corollary \ref{relation-toric}. Finally, in the last section, we discuss potential generalizations of these characterizations for $n \geq 3$.

Finally, in Section 5, we analyze the numerical semigroups that have extreme properties with respect to the inclusions exhibited in Proposition 2.2 and propose some conjectures supported by our results in previous sections and by computational experiments.

While working on this project, we have been making intensive and constant use of the \texttt{GAP} \cite{GAP} package \texttt{numericalsgps} \cite{numericalsgps} to generate examples and support conjectures. These examples allowed us to give shape to some of the results presented in this manuscript. New functionalities were added to \texttt{numericalsgps} to ease the development of this work, which can be consulted on the changelog of the package (changes from version 1.2.1 to 1.3.1). We also used the software \texttt{SageMath} \cite{sage} for computing Gr\"obner fans and universal Gr\"obner bases of polynomial ideals.

\section{Toric bases}\label{Section 2}

There are different sets of pure difference binomials associated to a toric ideal that contain relevant information about it, these are known as toric bases. In this section we define these sets and collect some of their properties for toric ideals correspoding to numerical semigroups. For a deeper and more general treatment of toric bases, see \cite{PET2, RTT, ST, TT, TTbases}.

Let $A = \{a_1, \ldots, a_n\} \subseteq \mathbb{Z}^+$. A binomial $\textbf{x}^{\textbf{u}}-\textbf{x}^{\textbf{v}}$ in $I_A$ is called \emph{primitive} if there is no other binomial $\textbf{x}^{\textbf{w}}- \textbf{x}^{\textbf{z}}$ in $I_A$, such that $\textbf{x}^{\textbf{w}}$ divides $\textbf{x}^{\textbf{u}}$ and $\textbf{x}^{\textbf{z}}$ divides $ \textbf{x}^{\textbf{v}}$. The set of primitive binomials, which is finite, is called the \emph{Graver basis} of $I_A$ and is denoted by $\operatorname{Gr}_A$. The \emph{universal Gr\"{o}bner basis} of $I_A$, denoted by $\mathcal{U}_A$, is the union of all reduced Gr\"obner bases of $I_A$; this set is finite and consists of binomials (see, for example, \cite{ST}). It is obviously a Gr\"obner basis of $I_A$ with respect to all monomial orders on $\mathbb{K}[\mathbf x]$. A \emph{Markov basis} of $I_A$ is a  binomial generating set of $I_A$ which is minimal for inclusion (its name comes from its relation with some Markov chains, see \cite[Theorem 3.1]{DST}). The \emph{universal Markov basis} of $I_A$, is denoted by $\mathcal{M}_A$, is the union of all the Markov bases of $I_A$, identifying the binomials with opposite signs. The elements of $\mathcal M_A$ are called \emph{minimal binomials}. Nakayama's lemma guarantees that all Markov bases have the same cardinality and that the $A$-degrees appearing in any Markov basis are invariant, these values are called \emph{Betti degrees} of $I_A$. We denote by $\betti{A}$ the set of the Betti degrees of $I_A$.

For $ i\in\{1,\dots,n\}$, we set 
\begin{equation}\label{eq:cb}
c_i := \min\left\{c \in \mathbb Z^+ \mid c\, a_i \in \langle a_1, \ldots, a_{i-1}, a_{i+1}, \ldots, a_n \rangle \right\}.
\end{equation}
If $c_i a_i = \sum_{j\neq i} r_{ij} a_j$ with $r_{ij} \in \N, j \in \{1,\ldots,n\} \setminus \{i\}$, then the binomial \[x_i^{c_i} - \prod_{j \in \{1,\ldots,n\} \setminus \{i\}}  x_j^{r_{ij}},\] is said to be a \emph{critical binomial of $I_A$ with respect to $x_i$}. The set of all critical binomials of $I_A$ is denoted by $\operatorname{Cr}_A$. The concept of critical binomial was introduced by Eliahou \cite{E} and later studied in \cite{AV} and \cite{KO}, among others. We emphasize here that $\{c_1 a_1, \ldots, c_n a_n\} \subseteq \betti{A}$ (see \cite[Proposition 2.3]{KO} for further details).

The \emph{support} of $\textbf{u} = (u_1,\ldots,u_n)\in \mathbb{N}^n$ is defined as $\operatorname{supp}(\textbf{u}):=\{ i \mid u_i\neq 0\}$. A binomial $\textbf{x}^{\textbf{u}}- \textbf{x}^{\textbf{v}} \in I_A$ is said to be {\it irreducible} if  $\operatorname{supp}(\textbf{u}) \cap \operatorname{supp}(\textbf{v}) = \emptyset$ and the nonzero entries of $\textbf{u} + \textbf{v}$ are relatively prime.
An irreducible binomial $\textbf{x}^{\textbf{u}}- \textbf{x}^{\textbf{v}} \in I_A$ is called a \emph{circuit} of $I_A$ if it has minimal support with respect to set inclusion, that is, there is not another binomial $\textbf{x}^{\textbf{u}'}- \textbf{x}^{\textbf{v}'} \in I_A$ with $\operatorname{supp}(\textbf{u}') \cup \operatorname{supp}(\textbf{v}') \subsetneq  \operatorname{supp}(\textbf{u}) \cup \operatorname{supp}(\textbf{v})$. The set of the circuits of $I_A$, denoted $\mathcal C_A$, can be explicitly described as follows (see, for example, \cite[Chapter 4]{ST} or \cite[Lemma 2.8]{KO}).

\begin{lem1}\label{lm:circuitsnumerical}
If $A = \{a_1, \ldots, a_n\}$ is a set of positive integers, then 
\[ \mathcal C_{A} = \left\{x_i^{a_j/\gcd(a_i,a_j)}  -  x_j^{a_i/\gcd(a_i,a_j)} \mid i,j\in \{1,\dots, n\},\ i \neq j \right\}, \]
\end{lem1}

The Graver basis, the universal Gr\"obner basis, the universal Markov basis, the set of the circuits and the set of the critical binomials are usually called \emph{toric bases}. The following inclusions of the previous sets are fulfilled.

\begin{prop1} \label{pr:inclusions} \cite[Proposition 4.11]{ST}, \cite[Theorem 2.3]{THO1} \cite[Proposition~2.3]{KO} If $A$ is a finite set of positive integers, then $$\mathcal C_A \subseteq \mathcal U_A \subseteq \operatorname{Gr}_A {\text \ and  \ } \operatorname{Cr}_A \subseteq \mathcal M_A \subseteq \operatorname{Gr}_A.$$
\end{prop1}

Some famous classes of ideals are defined through equality between some toric bases. These classes satisfy interesting geometric, combinatorial and homological properties. Let us show some of these classes with special emphasis on the case where $A$ is a finite subset of $\mathbb{Z}^+$. 

If every reduced Gr\"obner basis of $I_A$ consists of squarefree binomials, then $I_A$ is said to be \emph{unimodular} (see \cite[Remark 8.10]{ST}); in this case, $\mathcal C_A =  \operatorname{Gr}_A$. Note that if $I_A$ is unimodular and $A = \{a_1, \ldots, a_n\} \subseteq \mathbb{Z}^+$, then $\mathcal C_A \subseteq 
\mathcal{U}_A$ implies $\mathcal C_A = \{x_i - x_j \mid i \neq j\} \subseteq \mathcal{U}_A$, that is, $a_1 = a_2 = \ldots = a_n$. If $\mathcal U_A$ is a Markov basis of $I_A$, then $I_A$ is called \emph{robust} (see \cite{BR}). One of the most important properties of robust toric ideals is that the minimum number of generators of $I_A$ coincides with that of any of its initial ideals:
\begin{equation} \label{eq:muin} \mu(I_A) = \mu\big(\operatorname{in}_\prec (I_A)\big) \end{equation} for every monomial order $\prec$ on $\mathbb{K}[\mathbf{x}]$.
In \cite{T16} \emph{generalized robust} toric ideals are introduced. A toric ideal is generalized robust if $\mathcal{U}_A = \mathcal{M}_A$. In \cite[Theorem 4.12]{GT}, generalized robust toric ideals corresponding to numerical semigroups are characterized as those with a unique Betti degree. Finally, when $\operatorname{Gr}_A$ is a subset of every Markov basis of $I_A$, then $I_A$ is called \emph{strongly robust} (see \cite{Sull}). Observe that, in this case, one has that $\operatorname{Gr}_A = \mathcal{M}_A$.

\begin{figure}[ht]
\begin{center}
\begin{tikzpicture}[descr/.style={fill=white,inner sep=1.5pt}]
    \matrix (m) [matrix of math nodes,row sep=1em,column sep=7.1em]
  {
       \mathcal{C}_A & \operatorname{Gr}_A & \mathcal{M}_A & \mathcal{U}_A &   M_A 	\\		
		};
     \path
 (m-1-3) edge [double equal sign distance] node[sloped, above] {\small{generalized robust}}(m-1-4) 
						 (m-1-2)  edge [double equal sign distance] node[sloped, above] {\small{strongly robust}} (m-1-3)
						 (m-1-1)edge [double equal sign distance] node[sloped, above] {\small{unimodular}} (m-1-2)
						 (m-1-4) edge [double equal sign distance] node[sloped, above] {\small{robust}} (m-1-5);
\end{tikzpicture}
\caption{Some well-known classes of toric ideals, here $M_A$ denotes a Markov basis of $I_A$.}
\label{Fig.well-known-classes}
\end{center}
\end{figure}
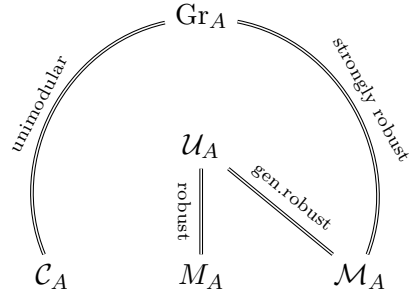

We end this section by summarizing how toric bases behave with respect to elimination of variables (see \cite[Proposition 4.13]{ST} for further details).
 
\begin{prop1}\label{pr:elimination} 
Let $A$ be a finite set of positive integers. If $A' \subseteq A$ is not empty, then
\begin{itemize}
\item[(a)] $I_{A'} = I_A \cap \mathbb K[\x_{A'}]$,
\item[(b)] $\mathcal C_{A'} = \mathcal C_A \cap \mathbb K[\x_{A'}]$,
\item[(c)] $\mathcal U_{A'} = \mathcal U_A \cap \mathbb K[\x_{A'}]$,
\item[(d)] $\operatorname{Gr}_{A'} = \operatorname{Gr}_A \cap \mathbb K[\x_{A'}]$,
\item[(e)] $\mathcal M_{A'} \supseteq \mathcal M_A \cap \mathbb K[\x_{A'}]$,
\item[(f)] $\operatorname{Cr}_{A'} \supseteq \operatorname{Cr}_A \cap \mathbb K[\x_{A'}]$,
\end{itemize}
where $\mathbb K[\x_{A'}] :=  \mathbb K[x_i \mid a_i \in A']$.
\end{prop1}

In general, we do not have equality in Proposition~\ref{pr:elimination}(e) and (f) as shown in the following example.

\begin{ex1}\label{notwell-elimin.ofMarkov-critical}
For $A' = \{a_1,a_2\} \subseteq A = \{a_1,a_2,a_3\} \subseteq \N$ with $a_1 = 4, a_2 = 5$ and $a_3 = 6$, we have that $\mathcal M_A = \operatorname{Cr}_A = \{x_1^3 - x_3^2, x_2^2 - x_1 x_3\}$. Thus, $\mathcal M_A \cap \mathbb K[x_1,x_2] =\operatorname{Cr}_A \cap \mathbb K[x_1,x_2]=\emptyset$, whereas  $\mathcal M_{A'} = \operatorname{Cr}_{A'} = \{x_1^5 - x_2^4\}$.

We can calculate the exponents of the binomials with the help of the \texttt{GAP} \cite{GAP} package \texttt{numericalsgps} \cite{numericalsgps}.
\begin{verbatim}
gap> AllMinimalRelationsOfNumericalSemigroup(NumericalSemigroup(4,5,6));
[ [ [ 1, 0, 1 ], [ 0, 2, 0 ] ], [ [ 3, 0, 0 ], [ 0, 0, 2 ] ] ]
gap> PrimitiveRelationsOfKernelCongruence([[4],[5],[6]]);
[ [ [ 0, 0, 2 ], [ 3, 0, 0 ] ], [ [ 0, 0, 3 ], [ 2, 2, 0 ] ], 
  [ [ 0, 0, 4 ], [ 1, 4, 0 ] ], [ [ 0, 0, 5 ], [ 0, 6, 0 ] ], 
  [ [ 0, 2, 0 ], [ 1, 0, 1 ] ], [ [ 0, 2, 1 ], [ 4, 0, 0 ] ],
  [ [ 0, 4, 0 ], [ 5, 0, 0 ] ] ]
\end{verbatim}
\begin{verbatim}
gap> AllMinimalRelationsOfNumericalSemigroup(NumericalSemigroup(4,5));
[ [ [ 5, 0 ], [ 0, 4 ] ] ]
gap> PrimitiveRelationsOfKernelCongruence([[4],[5]]);
[ [ [ 0, 4 ], [ 5, 0 ] ] ]
\end{verbatim}
\end{ex1}

As an easy consequence of Propositions \ref{pr:inclusions} and \ref{pr:elimination}, we obtain the following.

\begin{prop1}\label{genrobustfewvariables} 
Let $A$ be a finite set of positive integers and let $A'$ be a non-empty subset of $A$. 
\begin{itemize}
\item[(a)] If $\mathcal C_A \subseteq \mathcal M_A$, then $\mathcal C_{A'} \subseteq \mathcal M_{A'}$.
\item[(b)] If $\operatorname{Gr}_A =  \mathcal M_A$, then $\operatorname{Gr}_{A'} = \mathcal M_{A'}$.
\end{itemize} 
\end{prop1}




\section{Free and universally free numerical semigroups} 

Let $A = \{a_1, \ldots, a_n\}$ be a set of relatively prime positive integers and let $\mS$ be the numerical semigroup generated by $A$. When no proper subset of $A$ generates $\mS$, the set $A$ is said to be the \emph{minimal generating set} of $\mS$. Minimal generating sets of numerical semigroups always exist and are finite and unique (see \cite[Theorem 2.7]{RG-ns}). The cardinality of $A$ is known as the \emph{embedding dimension} of $\mS$.

Given a non-trivial partition $A_1$, $A_2$ of a subset of positive integers $A$ (that is, $A_1\cup A_2=A$, $A_1\cap A_2=\emptyset$ and both $A_1$ and $A_2$ are nonempty), we say that $A$ is the \emph{gluing} of $A_1$ and $A_2$, provided that $\lcm\big(\gcd(A_1),\gcd(A_2)\big)\in \langle A_1\rangle \cup \langle A_2\rangle$ (see \cite[Chapter~8]{RG-ns} for other equivalent definitions).

If $\mS$ is minimally generated by $A$, and $A$ is a gluing of $A_1$ and $A_2$, then we can write $\mS = d_1 \mS_1+ d_2 \mS_2$, where $d_1=\gcd(A_1)$, $d_2=\gcd(A_2)$, $\mS_1=\langle A_1/d_1\rangle$ and $\mS_2=\langle A_2/d_2\rangle$. In this setting, \cite[Theorem~10]{GO} asserts that
\begin{equation}\label{eq:betti-gluing}
    \betti{A}=\big\{d_1 b_1 \mid b_1\in \operatorname{Betti}(A_1/d_1)\big\}\cup \big\{d_2 b_2 \mid b_2\in \operatorname{Betti}(A_2/d_2)\big\} \cup \big\{d_1d_2\big\}.
\end{equation}
Observe that as $\gcd(A)=\gcd(d_1,d_2)=1$, we have that $d_1d_2=\lcm(d_1,d_2)$. 

Notice that if $A=\{a_1,\dots,a_n\}$, $A_1=\{a_1\}$ and $A_2=\{a_2,\dots,a_n\}$, then $A$ is a gluing of $A_1$ and $A_2$ if and only if $\lcm\big(a_1,\gcd(A_2)\big)\in \langle A_2\rangle$. In particular, as $\operatorname{Betti}(1)=\emptyset$, \eqref{eq:betti-gluing} translates in this setting to
\begin{equation}\label{eq:betti-gluing-one}
    \betti{A}=\big\{\gcd(A_2) b \mid b\in \operatorname{Betti}\big(A_2/\gcd(A_2)\big)\big\} \cup \big\{\lcm\big(a_1,\gcd(A_2)\big)\big\}.
\end{equation}

If $A$ is a gluing of $A_1$ and $A_2$, by using Delorme's correspondence \cite{D} and \cite[Theorem~9.2]{RG-ns} we can recover a minimal generating set of $I_A$ from those of $I_{A_1}$ and $I_{A_2}$. Next, we explicitly write this description for the case that $A_1$ is a singleton. This kind of gluing has been called \emph{monoscopic gluing} in \cite{mon-gl}.

\begin{prop1}\label{pr:prefree}
Let $A = \{ a_1,\ldots,a_n \}$ be a set of relatively prime positive integers, and set $d := \gcd(a_2,\ldots,a_n)$. If $da_1 = \alpha_2 a_2 + \cdots + \alpha_n a_n$ for some $\alpha_2,\ldots,\alpha_n~\in~\mathbb N$, then 
\[ I_A = I_{A \setminus \{a_1\}} \cdot \mathbb K[\x] + \langle x_1^{d} - x_2^{\alpha_2} \cdots x_n^{\alpha_n} \rangle.\]
\end{prop1}

\subsection{Free numerical semigroups}
Free numerical semigroups are those numerical semigroups obtained by gluing copies of $\mathbb{N}$.

\begin{def1}\label{freesemigroups}
A numerical semigroup $\mS$ with minimal generating set $A = \{a_1,\ldots,a_n\}$ is said to be \emph{free for the arrangement} $(a_1,\ldots,a_n)$ if for every $i\in \{1,\dots,n-1\}$, the set $\{a_i,\dots,a_n\}$ is the gluing of $\{a_i\}$ and $\{a_{i+1},\dots,a_n\}$, or equivalently, 
\begin{equation}\label{eq:free}  \lcm\big(a_i,\gcd(a_{i+1},\ldots,a_n)\big) \in \langle a_{i+1},\ldots, a_n\rangle, \text{ for all } i \in \{1,\ldots,n-1\}. 
\end{equation}  
\end{def1}

Notice that condition \eqref{eq:free} holds if and only if there exist $\alpha_{(i,i+1)},\ldots,\alpha_{(i,n)}\in\mathbb{N}$ such that  $\lcm\big(a_i,\gcd(a_{i+1},\ldots,a_n)\big) = \alpha_{(i,i+1)}a_{i+1} + \cdots + \alpha_{(i,n)}a_{n}$ for all $i \in \{1,\ldots,n-1\}$. 


\begin{cor1}\label{cor:betti-free}
Let $\mS$ be a numerical semigroup with minimal generating set $A=\{a_1,\dots, a_n\}$. If $\mS$ is free for the arrangement $(a_1,\dots,a_n)$, then \[\betti{A} = \big\{ \lcm\big(a_i,\gcd(a_{i+1}, \ldots, a_n)\big) \mid i\in \{1, \dots, n-1\}\big\}.\]
In particular, $\lcm(a_{n-1},a_n)\in \betti{A}$.
\end{cor1}
\begin{proof}
For $n=1$, we have that $\mS$ is $\mathbb{N}$ and thus it has no Betti degrees. For $n=2$, we apply \eqref{eq:betti-gluing-one} obtaining $\betti{A}=\emptyset \cup \{\lcm(a_1,a_2)\}$. Now assume that the result holds for numerical semigroups with embedding dimension $n-1$, and let us show that it holds for embedding dimension $n$. So assume that $\mS$ is free for the arrangement $(a_1,\dots,a_n)$ of its minimal generators. Let $\mS'$ be the numerical semigroup minimally generated by the set $A'=\{a_2/d,\dots,a_n/d\}$, with $d=\gcd(a_2,\dots,a_n)$. By definition, $\mS'$ is also free for the arrangement $(a_2/d,\dots,a_n/d)$ of its minimal generators. Again, considering that $A$ is a gluing of $\{a_1\}$ and $\{a_2,\ldots,a_n\}$, by using \eqref{eq:betti-gluing-one}, we obtain that $$\betti{A}=\big\{ d b \mid b\in  \betti{A'}\big\}\cup \big\{\lcm\big(a_1, \gcd(A')\big)\big\}.$$ By induction hypothesis $$\betti{A'}=\big\{ \lcm\big(a_i/d,\gcd(a_{i+1}/d, \ldots, a_n/d)\big) \mid i\in \{2, \dots, n-1\}\big\},$$ and so $$\big\{ d b \mid b\in  \betti{A'}\big\}= \big\{ \lcm\big(a_i,\gcd(a_{i+1}, \ldots, a_n)\big) \mid i\in \{2, \dots, n-1\}\big\},$$ which concludes the proof.

The last assertion is just the case $i=n-1$.
\end{proof}



Recall that a numerical semigroup is said to be \emph{free} if it is free for an arrangement of its minimal generating set (see \cite[Section 2.3]{AG}). 

The following result is just a reformulation of the description of a minimal presentation for a free numerical semigroup given 
\cite[Corollary~9.19]{RG-ns}, and can be obtained by applying Proposition~\ref{pr:prefree} inductively. We warn the reader that the arrangement of the generators in 
\cite{RG-ns} is taken in the reverse order to that taken by us.

\begin{thm1}\label{th:free}
Let $\mS$ be the numerical semigroup minimally generated by $A = \{a_1,\ldots,a_n\}$. Then $\mathcal S$ is free for the arrangement $(a_1, \ldots, a_n)$ if and only if there exist positive integers $r_i$, for every $i \in \{1, \dots, n-1\}$, and non-negative integers $r_{ij}$, for every $j\in \{ i+1, \ldots, n\}$ and $i\in\{ 1, \ldots, n-1\}$, such that the set \[\left\{x_i^{r_i}-\prod_{j=i+1}^n x_j^{r_{ij}} \mid i \in\{ 1, \dots, n-1\}\right\}\] is a Markov basis of $I_A$. 
\end{thm1}

As a direct consequence of Theorem~\ref{th:free}, every free numerical semigroup is a complete intersection (see also \cite[Corollary~9.18]{RG-ns}). When $n = 2$ every numerical semigroup is free. For $n=3$, Herzog proved in \cite{H} that $\mS$ is free if and only if $I_A$ is a complete intersection. For $n \geq 4$, there are complete intersection numerical semigroups which are not free as it is shown in the following example.

\begin{ex1}
The numerical semigroup $\mS = \langle 10,14,15,21\rangle$ is not free for any arrangement of the generators. On the other hand, the ideal $I_A$ is a complete intersection; indeed, one can check that $\{ x_1^3 - x_3^2, x_2^3 - x_4^2, x_1^2x_3 - x_2x_4\}$ is a Markov basis for $I_A$ (and the universal Markov basis of $I_A$). 
\begin{verbatim}
gap> s:=NumericalSemigroup(10,14,15,21);;
gap> IsCompleteIntersection(s);
true
gap> AsGluingOfNumericalSemigroups(s);
[ [ [ 10, 15 ], [ 14, 21 ] ] ]
gap> IsFree(s);
false
gap> MinimalPresentation(s);
[ [ [ 0, 0, 0, 2 ], [ 0, 3, 0, 0 ] ], 
  [ [ 0, 0, 2, 0 ], [ 3, 0, 0, 0 ] ],
  [ [ 0, 1, 0, 1 ], [ 2, 0, 1, 0 ] ] ]
 \end{verbatim}
\end{ex1}

\subsection{Universally free numerical semigroups}

Next, we recall the notion of a universally free numerical semigroup which we study in the rest of the paper.

\begin{def1}
A numerical semigroup $\mS$ is called \emph{universally free}, if it is free for any arrangement of its minimal generating set.
\end{def1}

Observe that free numerical semigroups are not necessarily universally free as it is shown in the following example.

\begin{ex1}\label{ex:free}
Consider the numerical semigroup $\mS = \langle a_1,a_2,a_3,a_4 \rangle$ with $a_1 = 8,\, a_2 = 9,\, a_3 = 10,\, a_4 = 12$. We have that $\mS$ is not free for the arrangement $(a_1,a_2,a_3,a_4)$ because $\lcm\big(a_1,\gcd(a_2,a_3,a_4)\big) = 8 \notin \langle a_2,a_3,a_4\rangle$. Hence, it is not universally free . However, $\mS$ is free for the arrangement $(a_2, a_3, a_1, a_4)$. Indeed, \begin{itemize} 
\item $\lcm\big(a_2,\gcd(a_3,a_1,a_4)\big) = 18 = 2a_2= a_1 + a_3 \in \langle a_1,a_3,a_4\rangle$,
\item $\lcm\big(a_3,\gcd(a_1,a_4)\big) = 20 = 2a_3 = a_1 + a_4 \in \langle a_1,a_4 \rangle$, and
\item $\lcm(a_1,a_4) = 24 = 3a_1 = 2 a_4 \in \langle a_4 \rangle$.
\end{itemize}
Thus, $\mS$ is a free numerical semigroup. 

It follows that a Markov basis for $I_A$ is $\{x_2^2 - x_1x_3,\, x_3^2 - x_1x_4, x_1^3-x_4^2\}.$ We remark that the minimal generator $x_3^2 - x_1x_4$ is not a circuit. Here, $I_A$ is not generated by its circuits. As we will prove in the next proposition, this does not happen in the case of universally free numerical semigroups.
\end{ex1}



\begin{rem1}
Universally free numerical semigroups are not so common among free semigroups, not even among telescopic numerical semigroups. Recall that a numerical semigroup minimally generated by $\{a_1,\dots,a_n\}$ is said to be telescopic if it is free for the arrangement $(a_n,\ldots,a_1)$ where $a_n>\dots >a_1$ (see for example \cite[Section~2.3]{AG}); both terms free and telescopic were originally the same, and discovered independently \cite{bertin,KirPel}, but at some point some authors used the term telescopic for freeness with respect to the arrangement given by the natural order of the generators.

Below we give a table with the number of free ($nf_f$), telescopic ($nt_f$) and universally free ($nsf_f$) numerical semigroups with Frobenius number $f\in\{101,111,\dots,691\}$.
\begin{center}
\begin{tabular}{c|ccc}
$f$ & $nf_f$ & $nt_f$ & $nsf_f$ \\ \hline
101 & 194 & 86 & 5\\ 
111 & 169 & 83 & 3\\ 
121 & 310 & 140 & 3\\ 
131 & 387 & 171 & 7\\ 
141 & 330 & 151 & 2\\ 
151 & 571 & 230 & 6\\ 
161 & 667 & 281 & 5\\ 
171 & 588 & 257 & 2\\ 
181 & 949 & 367 & 5\\ 
191 & 1130 & 414 & 8\\ 
201 & 938 & 395 & 2\\ 
211 & 1502 & 565 & 7\\ 
221 & 1742 & 639 & 7\\ 
231 & 1444 & 569 & 2\\ 
241 & 2284 & 844 & 7\\ 
251 & 2602 & 928 & 12\\ 
261 & 2194 & 831 & 2\\ 
271 & 3337 & 1120 & 7\\ 
281 & 3748 & 1293 & 5\\ 
291 & 3113 & 1140 & 3\\ 
\end{tabular}
\begin{tabular}{c|ccc}
$f$ & $nf_f$ & $nt_f$ & $nsf_f$ \\ \hline
301 & 4682 & 1538 & 3\\
311 & 5266 & 1697 & 10\\
321 & 4361 & 1568 & 4\\
331 & 6515 & 2099 & 6\\
341 & 7205 & 2279 & 8\\
351 & 6038 & 2035 & 3\\
361 & 8781 & 2789 & 5\\
371 & 9629 & 3016 & 9\\
381 & 7964 & 2638 & 4\\
391 & 11631 & 3507 & 8\\
401 & 12763 & 3928 & 10\\
411 & 10519 & 3408 & 2\\
421 & 15165 & 4504 & 5\\
431 & 16504 & 4798 & 16\\
441 & 13738 & 4425 & 1\\
451 & 19503 & 5688 & 6\\
461 & 21105 & 6074 & 11\\
471 & 17412 & 5306 & 4\\
481 & 24744 & 7117 & 7\\
491 & 26726 & 7557 & 12\\
\end{tabular}
\begin{tabular}{c|ccc}
$f$ & $nf_f$ & $nt_f$ & $nsf_f$ \\ \hline
501 & 22021 & 6618 & 3\\
511 & 30903 & 8495 & 7\\
521 & 33383 & 9291 & 11\\
531 & 27895 & 8240 & 5\\ 
541 & 38575 & 10443 & 5\\ 
551 & 41361 & 11052 & 12\\ 
561 & 34063 & 9986 & 2\\ 
571 & 47422 & 12716 & 10\\ 
581 & 50510 & 13445 & 6\\ 
591 & 41842 & 11732 & 5\\ 
601 & 57902 & 15424 & 10\\ 
611 & 61733 & 16207 & 14\\ 
621 & 51545 & 14350 & 2\\ 
631 & 70050 & 17974 & 13\\ 
641 & 74552 & 19425 & 11\\ 
651 & 61350 & 16919 & 3\\ 
661 & 84303 & 21419 & 9\\ 
671 & 89396 & 22469 & 11\\ 
681 & 73826 & 20148 & 3\\ 
691 & 100705 & 25393 & 9\\ 
\end{tabular}
\end{center}
This table has been constructed by using \texttt{TelescopicNumericalSemigroupsWith\-Frobenius\-Number} and \texttt{FreeNumericalSemigroupsWithFrobeniusNumber} from the \texttt{GAP} \cite{GAP} package \texttt{numerical\-sgps} \cite{numericalsgps}. The number of universally free numerical semigroups has been obtained by filtering those telescopic numerical semigroups that are universally free.
\end{rem1}

In the next proposition, we prove some very useful (for the rest of the present manuscript) properties of a universally free numerical semigroup.

\begin{prop1}\label{prop:free}
Let $\mS$ be the numerical semigroup with minimal generating set $A = \{a_1,\ldots,a_n\}$, and let $c_i$ be as in \eqref{eq:cb}. If $\mathcal S$ is universally free, then
\begin{itemize}
    \item[(a)] $c_i = \gcd(A \setminus \{a_i\}),$ for every $i \in \{1, \dots, n\}$;
    \item[(b)] each Betti degree of $I_A$ is divisible by $c_k a_k$ for some $k \in \{1, \ldots, n\}$;
    \item[(c)] $I_A$ is minimally generated by $n-1$ circuits;
    \item[(d)] $\betti{A} = \{\lcm(a_i, a_j) \mid i \neq j\}$.
\end{itemize}
\end{prop1}
\begin{proof}
(a) Given $A=\{a_1,\ldots,a_n\}\subseteq \mathbb{Z}^+$, let us denote by $\tilde{c}_i$ the minimum positive integer $k$ such that $ka_i$ is in the group spanned by $A\setminus\{a_i\}$, that is:
\begin{equation}\label{eq:tilde-c}
\tilde{c}_i=\min\left\{ k\in \mathbb{Z}^+ \mid k a_i\in \sum_{j\in \{1,\dots,n\}\setminus\{i\}} \mathbb{Z}a_j\right\}.     
\end{equation} We claim that $\tilde{c}_i=c_i$ for all $i=1,\ldots,n$.

Notice that $$\langle A\setminus\{a_i\}\rangle =  \sum_{j\in \{1,\dots,n\}\setminus\{i\}}\mathbb{N} a_j \subseteq \sum_{j\in \{1,\dots,n\}\setminus\{i\}} \mathbb{Z} a_j,$$ and then it follows that $\tilde{c}_i\le c_i$. 

Also, the group spanned by $A\setminus\{a_i\}$ is precisely the group spanned by $\gcd(A\setminus\{a_i\})$, and so $\tilde{c}_i$ is the minimum positive integer $k$ such that $ka_i$ is a multiple of $\gcd(A\setminus\{a_i\})$. Thus, $\tilde{c}_i a_i$ is precisely $\lcm\big(a_i,\gcd(A\setminus\{a_i\})\big)$. Since $\mS$ is a universally free numerical semigroup minimally generated by $A$, by \eqref{eq:free}, we have that $\lcm\big(a_i,\gcd(A\setminus\{a_i\})\big)=\tilde{c}_i a_i\in \langle A\setminus \{a_i\}\rangle$, and in particular, $c_i\le \tilde{c}_i$. Thus, for universally free numerical semigroups $\tilde{c}_i=c_i$ for all $i=1,\ldots,n$. 

By the definition of a numerical semigroup we have that $\gcd(A)=1$ and so $$\lcm\big(a_i,\gcd(A\setminus\{a_i\})\big)=a_i\gcd(A\setminus\{a_i\}).$$ It follows that $c_i=\gcd(A\setminus\{a_i\}).$

(b) Let $\beta \in \betti{A}$. On the one hand, by Corollary~\ref{cor:betti-free}, $\beta = \lcm\big(a_k,\gcd(a_{k+1}, \ldots, a_n)\big)$ for some $k \in \{1, \ldots, n\}$; in particular, $\beta = d_k a_k$ for some $k\in\{1, \ldots, n\}$. On the other hand, since $\mS$ is universally free, it is free for the arrangement $(a_{\sigma(1)}, \ldots, a_{\sigma(n)})$ of $A$ where $\sigma$ is a permutation of $\{1, \ldots, n\}$ such that $\sigma(1) = k$. Now, arguing as in part (a), we conclude that $d_k a_k$ is divisible by $c_k a_k = \lcm\big(a_k,\gcd(A\setminus \{a_k\})\big)$.

(c) Let us use induction on $n$. For $n=2$, the result holds trivially, since the only Betti degree is $a_1a_2$, and the only generator of $I_A$ is $x_1^{a_2}-x_2^{a_1}$. 

Suppose that the statement is true for all universally free numerical semigroups with embedding dimension $n-1$ and let us show that it holds for $\mS$ minimally generated by the set $A=\{a_1,\dots,a_n\}$. We already know that $c_i a_i\in \betti{A}$ (see Proposition~\ref{pr:inclusions}), and that $\mS$ is a complete intersection (since it is free), and thus there must be $i\neq j$ such that $c_ia_i=c_ja_j$ (every Betti degree appears as the $A$-degree of at least one minimal generator of $I_A$, and this ideal has $n-1$ minimal generators). After rearranging the generators, we may assume that $c_1a_1=c_2a_2$. Let $A'=\{a_2/c_1,\dots,a_n/c_1\}$, recall that 
$c_1=\gcd(a_2,\dots,a_n)$ by (a). As $\mS$ is free for the arrangement $(a_1,a_2,\dots,a_n)$, by Proposition~\ref{pr:prefree}, a Markov basis of $I_A$ is obtained by adjoining the binomial $x_1^{c_1}-x_2^{c_2}$ to a minimal generating set of $I_{A'}$. 

By induction hypothesis we have that $I_A'$ is minimally generated by $n-2$ circuits. Let us prove that $x_1^{c_1}-x_2^{c_2}$ is a circuit. Denote by $m_{12}=\lcm(a_1,a_2)$ and by $d_{12}=\gcd(a_1,a_2)$. Since $(a_2/d_{12})a_1=(a_1/d_{12})a_2$, we have, by the minimality of $c_1$, that $c_1\le a_2/d_{12}$. Also, as $c_1a_1=c_2a_2$, we obtain that $c_1a_1$ is a common multiple of $a_1$ and $a_2$, and consequently $c_1a_1\ge m_{12} =(a_2/d_{12})a_1$, which in particular means that $c_1\ge a_2/d_{12}$. Hence $c_1=a_2/d_{12}$, and thus $c_2
= a_1/d_{12}$. By Lemma \ref{lm:circuitsnumerical} this proves that $x_1^{c_1}-x_2^{c_2}$ is a circuit and that $I_A$ is (minimally) generated by $n-1$ circuits, those inherited from $I_{A'}$ plus $x_1^{c_1}-x_2^{c_2}$.

(d) Since $\mS$ is free for any arrangement, the inclusion $ \{ \lcm(a_i,a_j)\mid i\neq j \} \subseteq \betti{A}$ follows by applying Corollary~\ref{cor:betti-free} to the arrangements where $a_i$ and $a_j$ are the last two elements of the corresponding arrangement. Since the $A$-degree of the circuit $x_i^{a_j/\gcd(a_i,a_j)} - x_j^{a_i/\gcd(a_i,a_j)}$ is $\lcm(a_i,a_j)$, the other inclusion follows from (c). 
\end{proof}

Let us see now that the conditions stated in Proposition \ref{prop:free} are not sufficient to be universally free.

\begin{ex1} 
Let $\mS$ be the semigroup generated by $\{30, 105, 546, 770\}$. From
\begin{verbatim}
gap> s:=NumericalSemigroup( 30, 105, 546, 770 );
gap> AsGluingOfNumericalSemigroups(s);
[ [ [ 30 ], [ 105, 546, 770 ] ], [ [ 30, 105, 546 ], [ 770 ] ], 
  [ [ 30, 105, 770 ], [ 546 ] ], [ [ 30, 546 ], [ 105, 770 ] ], 
  [ [ 30, 546, 770 ], [ 105 ] ], [ [ 30, 770 ], [ 105, 546 ] ] ]  
\end{verbatim}
we have that $c_i=\gcd(A\setminus\{a_i\})$ for all $i\in\{1,\ldots,4\}$. So Condition (a) in Proposition~\ref{prop:free} holds for this semigroup. The same stands for Condition (b):
\begin{verbatim}
gap> A:=MinimalGenerators(s);
[ 30, 105, 546, 770 ]
gap> Set(A, a-> a*Gcd(Difference(A,[a])));
[ 210, 2310, 2730 ]
gap> BettiElements(s);
[ 210, 2310, 2730 ]    
\end{verbatim}
However, $\lcm(546,770) = 30030 \notin \betti{A}.$ 
\begin{verbatim}
gap>  Set(Combinations(A,2),Lcm);
[ 210, 2310, 2730, 30030 ]
\end{verbatim}
Hence, Condition (d) in Proposition~\ref{prop:free} does not hold for $\mS$, and $\mS$ is not universally free. As a consequence, Conditions (a) and (b) in Proposition~\ref{prop:free} are not sufficient.

\end{ex1}

\begin{ex1}
There are many examples of numerical semigroups for which Condition (d) in Proposition~\ref{prop:free} holds, and they are not universally free. A search with \texttt{numericalsgps} for semigroups generated by four integers ranging between 10 and 500 throws many examples. For instance, $\langle  30, 40, 45, 72\rangle$ is not universally free but Condition (d) in Proposition~\ref{prop:free} holds.
\end{ex1}

One can modify the previous example to find a numerical semigroup with embedding dimension five which is not universally free and  fulfills conditions (a) to (d).

\begin{ex1}Consider the numerical semigroup $\mS$ minimally generated by the set $A = \{a_1,a_2,a_3,a_4,a_5\}$ with $a_1 = 30 \cdot 7 = 210,\, a_2 = 30 \cdot 11 = 330,\, a_3 = 40 \cdot 7 \cdot 11 =  3080,\, a_4 = 45 \cdot 7 \cdot 11 = 3465,\, a_5 = 72 \cdot 7 \cdot 11 = 5544$. One can check that $\mS$ is not free for the arrangement $(a_1,a_2,a_3,a_4,a_5)$ since $\lcm\big(a_2,\gcd(a_3,a_4,a_5)\big) = 7 a_2 = 2310 \notin \langle a_3,a_4,a_4\rangle$ and, hence, it is not universally free. Nevertheless $\mS$ is free for the arrangement $(a_5,a_4,a_3,a_2,a_1)$; indeed,
\begin{itemize}
    \item $\lcm\big(a_5,\gcd(a_1,a_2,a_3,a_4)\big) = 5 a_5 = 8 a_4 \in \langle a_1,a_2,a_3,a_4\rangle,$
    \item $\lcm\big(a_4,\gcd(a_1,a_2,a_3)\big) = 2 a_4 = 21 a_2 \in \langle a_1,a_2,a_3 \rangle,$
    \item $\lcm\big(a_3,\gcd(a_1,a_2)\big) = 3 a_3 = 28 a_2 \in \langle a_1,a_2 \rangle,$ and
    \item $\lcm(a_2,a_1) = 7 a_2 = 11 a_1.$
\end{itemize}
Thus, a Markov basis for $I_A$ is \[ \{x_5^5 - x_4^8,\, x_4^2 - x_2^{21},\, x_3^{3} - x_2^{28},\, x_2^7 
- x_1^{11}\}\]
and $I_A$ is generated by $4$ circuits.
We have that $c_1 = 11,\, c_2 = 7,\, c_3 = 3,\,c_4 = 2\, c_5 = 5$, and 
\[ \betti{A} = \{ 2310, 6930, 9240, 27720 \} = \{c_2a_2, c_4a_4, c_3a_3, c_5a_5\} = \{ \lcm(a_i,a_j) \, \vert \, i \neq j\}. \]
\begin{verbatim}
gap> t:=NumericalSemigroup(30*7, 30*11, 40*77, 45*77, 72*77);
<Numerical semigroup with 5 generators>
gap> IsUniversallyFree(t);
false
gap> IsFree(t);
true
gap> AsGluingOfNumericalSemigroups(t);
[ [ [ 210 ], [ 330, 3080, 3465, 5544 ] ], 
  [ [ 210, 330, 3080 ], [ 3465, 5544 ] ],
  [ [ 210, 330, 3080, 3465 ], [ 5544 ] ], 
  [ [ 210, 330, 3080, 5544 ], [ 3465 ] ],
  [ [ 210, 330, 3465 ], [ 3080, 5544 ] ], 
  [ [ 210, 330, 3465, 5544 ], [ 3080 ] ],
  [ [ 210, 3080 ], [ 330, 3465, 5544 ] ], 
  [ [ 210, 3080, 3465 ], [ 330, 5544 ] ],
  [ [ 210, 3080, 3465, 5544 ], [ 330 ] ], 
  [ [ 210, 3080, 5544 ], [ 330, 3465 ] ],
  [ [ 210, 3465 ], [ 330, 3080, 5544 ] ], 
  [ [ 210, 3465, 5544 ], [ 330, 3080 ] ],
  [ [ 210, 5544 ], [ 330, 3080, 3465 ] ] ]
gap> A:=MinimalGenerators(t);
[ 210, 330, 3080, 3465, 5544 ]
gap> Set(Combinations(A,2),Lcm);
[ 2310, 6930, 9240, 27720 ]
gap> BettiElements(t);
[ 2310, 6930, 9240, 27720 ]
gap> MinimalPresentation(t);
  [ [ [ 0, 0, 0, 0, 5 ],[ 0, 0, 0, 8, 0 ] ], 
    [ [ 0, 0, 0, 2, 0 ],[ 0, 21, 0, 0, 0 ] ],
    [ [ 0, 0, 3, 0, 0 ],[ 0, 28, 0, 0, 0 ] ], 
    [ [ 0, 7, 0, 0, 0 ],[ 11, 0, 0, 0, 0 ] ] ]
\end{verbatim}
\end{ex1}

We are not aware of any numerical semigroup $\mS = \langle a_1,a_2,a_3,a_4\rangle$ which is not universally free and satisfies conditions (a), (b), (c) and (d). We wonder if these conditions are sufficient for being universally free in embedding dimension four.

We also want to highlight that, as the following example shows, there are universally free numerical semigroups where does not hold that every Betti degree of $I_A$ is of the form $c_ka_k$.
\begin{ex1}\label{ex:FAnotBettidiv0} 
Consider $\mS = \langle a_1 = 390, a_2 = 546, a_3 = 770, a_4 = 1155 \rangle$. Then one can check that $\mS$ is universally free. We have that $$I_A = \langle x_1^7 - x_2^5, x_3^3 - x_4^2, x_2^{55} - x_4^{26} \rangle,$$ and, hence, the Betti degrees are $\beta_1  = 2310 = c_3 a_3 = c_4 a_4,\ \beta_2 = 2730 = c_1 a_1 = c_2 a_2$ and $\beta_3 = 30030 \notin \{c_i a_i \, \vert \, 1 \leq i \leq 4\}$.
\begin{verbatim}
gap> s:=NumericalSemigroup(390,546,770,1155);;
gap> IsUniversallyFree(s);
true
gap> BettiElements(s);
[ 2310, 2730, 30030 ]
gap> A:=MinimalGenerators(s);;
gap> Set(A,g->Lcm(g,Gcd(Difference(A,[g]))));
[ 2310, 2730 ]
\end{verbatim}
Notice that, as we already know, the set of all Betti degrees coincides with the set of all $\lcm(a_i,a_j)$, $i\neq j$.
\begin{verbatim}
gap> Set(Combinations(A,2), Lcm);
[ 2310, 2730, 30030 ]
\end{verbatim}
\end{ex1}

\begin{prop1}\label{prop:free2}
Let $\mS$ be a numerical semigroup minimally generated by $A = \{a_1,\ldots,a_n\}$. 
If $\mS$ is universally free, 
then $\mathcal C_A \subseteq \mathcal M_A$.
\end{prop1}

\begin{proof}
Suppose that $\mS$ is universally free. 
Consider the circuit $$f = x_i^{a_j/\gcd(a_i,a_j)}  -  x_j^{a_i/\gcd(a_i,a_j)},\ i \neq j$$ and an arrangement $(a_{\sigma(1)}, \ldots, a_{\sigma(n)})$ of $A$ where $\sigma$ is a permutation of $\{1, \ldots, n\}$ such that $\sigma(n-1) = i$ and $\sigma(n) = j$. Since $S$ is universally free, it follows that it is free for the arrangement $(a_{\sigma(1)}, \ldots, a_{\sigma(n)})$. So, by Theorem \ref{th:free}, there exists a Markov basis of $I_A$ containing $x_i^{r_i} - x_j^{r_{ij}}$ for some $r_i > 0$ and $r_{ij} > 0$. Now, since $\gcd(r_i,r_{ij}) = 1$ and $\big(a_i/\gcd(a_i,a_j)\big) r_i = \big( a_j/\gcd(a_i,a_j)\big) r_{ij}$, we conclude that $r_{ij} = a_i/\gcd(a_i,a_j)$ and $r_i = a_j/\gcd(a_i,a_j)$. Therefore, $f \in \mathcal{M}_A$.
\end{proof}

Computational evidence supports that the opposite implication in the above proposition is also true (we had a process running for 20 days, checking a total of 65730378 combinations of four potential generators between 10 and 500 --yielding 48130628 numerical semigroups of embedding dimension four--, and whenever $\mathcal C_A \subseteq \mathcal M_A$, we got that $\langle A\rangle$ was a universally free numerical semigroup). 

\begin{conj1}\label{cj:conj1} Let $\mS$ be a numerical semigroup minimally generated by $A = \{a_1,\ldots,a_n\}$. Then, $\mS$ is universally free if and only if $\mathcal C_A \subseteq \mathcal M_A$.
\end{conj1}

Now we study the initial ideals of the toric ideal of a universally free numerical semigroup. The main result in this direction is Theorem \ref{universally-Grobner} where we prove that a numerical semigroup is universally free if and only if all the initial ideals of the corresponding toric ideal are complete intersections or, equivalently, all its minimal Gr\"obner bases have $n-1$ elements. To deduce this, we will use the following result from \cite{GT}.

\begin{prop1}\label{pr:freeGB}\cite[Proposition 4.5]{GT}  Let $\mS$ be the numerical semigroup with minimal generating set $A = \{a_1,\ldots,a_n\}$. Then $\mS$ is free for the arrangement $(a_1,\ldots,a_n)$ if and only if the reduced Gr\"obner basis of $I_A$ with respect to the lexicographic order on $\mathbb K[\mathbf x]$ with $x_1 \succ \cdots \succ x_n$ has $n-1$ elements.
\end{prop1}


\begin{thm1}\label{universally-Grobner}
Let $\mS$ be a numerical semigroup with minimal generating set $A = \{a_1,\ldots,a_n\}$. Then, $\mS$ is universally free if and only if every reduced Gr\"obner basis of $I_A$ has $n-1$ elements.
\end{thm1}
\begin{proof}
{\emph{Necessity.}} Suppose that $\mS$ is universally free and consider $\prec$ any monomial order. We are going to construct a (possibly not reduced) Gr\"obner basis of $I_A$  with respect to $\prec$ consisting of $n-1$ circuits. We know that $c_i a_i\in \betti{A}$ for all $i \in \{1,\ldots,n\}$ (see Proposition~\ref{pr:inclusions}), and that $\mS$ is a complete intersection (since it is free). Thus there must be $i\neq j$ such that $c_ia_i=c_ja_j$. We consider the binomial $f_1 := x_i^{c_i} - x_j^{c_j} \in I_A$ and assume without loss of generality that $x_i^{c_i} \succ x_j^{c_j}$ and that $i = 1$. Since $c_1 = d_1 = \gcd(a_2,\ldots,a_n)$, by Proposition \ref{pr:prefree} it follows that $I_A = \big(I_A \cap \mathbb K[x_2,\ldots,x_n]\big) \mathbb K [\bx] + (f_1)$. Moreover, $I_A \cap \mathbb K[x_2,\ldots,x_n]$ equals the toric ideal $I_{A'}$  being $A' = \{a_2/d_1,\ldots,a_n/d_1\}$. If we consider $\mS' = \langle A' \rangle$, then $\mS'$ is also universally free. Hence, one can iterate this procedure to get a set of generators  $\mathcal G = \{f_1,\ldots,f_{n-1}\}$ formed by circuits and whose leading terms are relatively prime (because they involve different variables). Thus,  $\mathcal G = \{f_1,\ldots,f_{n-1}\}$ is a Gr\"obner basis of $I_A$ with respect to $\prec$.

{\emph{Sufficiency.}} Assume that every reduced Gr\"obner basis of $I_A$ has $n-1$ elements. In particular every lexicographic Gr\"obner basis of $I_A$ has $n-1$ elements and, by Proposition \ref{pr:freeGB} it follows that $\mS$ is universally free. 
\end{proof}



In light of the last two results, it seems natural to ask the following question.

\begin{conj1} Let $\mS$ be a numerical semigroup minimally generated by $A = \{a_1,\ldots,a_n\}$. Then, $\mS$ is universally free if and only if $\mu(I_A) = \mu\big(\operatorname{in}_\prec(I_A)\big)$ for every monomial order $\prec$.
\end{conj1}

\subsection{Betti divisible numerical semigroups}\mbox{}\par

In \cite{GOR} the family of submonoids of $\N^n$ having unique Betti degree is studied. In particular, it is proved that numerical semigroups with a unique Betti degree are those in which the set of the circuits agrees with the set of the critical binomials and the Graver basis of the corresponding toric ideal. This implies that all inclusions in Proposition~\ref{pr:inclusions} are equalities. In particular, numerical semigroups having unique Betti degree are universally free.

The family of submonoids of $\N^n$ having unique Betti degree is generalized in \cite{GH}, where the authors study the family of Betti divisible submonoids of $\N^n$. An affine monoid $\mS \subseteq \N^n$ with minimal generating set $A$ is Betti divisible if its Betti degrees are ordered by divisibility. In the particular context of numerical semigroups, they characterize the sets of generators of numerical semigroups that are Betti divisible. In their description they implicitly use the following unique writing of any set of relatively prime positive integers.

\begin{lem1}\label{lm:uniquewriting} Every set $A = \{a_1,\ldots,a_n\}$ of relatively prime positive integers can be uniquely written as $a_i = f_i \prod_{j \neq i} d_j$  with the following conditions:
\begin{itemize} 
\item[(1)] $d_1,\ldots,d_n$ are pairwise coprime, 
\item[(2)] every subset of $n-1$ elements among $f_1,\ldots,f_n$ are relatively prime,
\item[(3)] $\gcd(f_i,d_i) = 1$ for all $i \in \{1,\ldots,n\}.$ 
\end{itemize}
\end{lem1}
\begin{proof}
 For the existence it suffices to take $d_i := \gcd(A-\{a_i\})$ and $f_i := a_i / \prod_{j\neq i} d_j$. To prove the uniqueness one just has to observe that if $a_i= f_i \prod_{j \neq i} d_j$ with $f_i$ and $d_i$ satisfying (1), (2) and (3), then $\gcd(A-\{a_i\}) = d_i$.
\end{proof}

\begin{prop1}\label{pr:shapeBettidiv} \cite[Theorem 7.10]{GH} 
A numerical semigroup $\mS$ is Betti divisible if and only if it is minimally generated by $A = \{a_1,\ldots,a_n\}$, where $a_i = f_i \prod_{j \neq i} d_j$ for all $i \in \{1,\ldots,n\}$, being $d_i, f_i$ some positive integers satisfying that
\begin{enumerate}
\item[(a)]  $d_1,\ldots,d_n$ are pairwise prime, 
\item[(b)] $1 = f_1 = f_2$ and $f_i$ divides $f_{i+1}$ for all $i \in \{2,\ldots,n-1\}$, 
\item[(c)] $\gcd(f_i,d_i) = 1$ for all $i \in \{1,\ldots,n\}$.
\end{enumerate}
\end{prop1}

In \cite[Theorem 7.12]{GH} it is stated that a numerical semigroup is Betti divisible if and only if it is free for every arrangement of its generators. However, there is a flaw in the proof and there are universally free semigroups that are not Betti divisible. 
One counterexample is the one provided in Example \ref{ex:FAnotBettidiv0}, in which the numerical semigroup is universally free, with Betti degrees $\{2310,\, 2730,\, 30030\}$ and therefore it is not Betti divisible.  

In spite of the above example, the implication a) implies b) in  \cite[Theorem 7.12]{GH} is true, that is to say, a Betti divisible numerical semigroup is free for every arrangement of its generators. Thus, Betti divisible numerical semigroups provide an interesting family of universally free numerical semigroups. 

\begin{thm1}\label{th:bettidivbases}If $\mS$ is a Betti divisible numerical semigroup minimally generated by $A$, then $\mathcal C_A = {\mathcal U}_A \subseteq \operatorname{Cr}_A = \mathcal M_A = \operatorname{Gr}_A$.
\end{thm1}
\begin{proof} Let $\mS$ be a Betti divisible numerical semigroup. Thus, its minimal set of generators $A = \{a_1,\ldots,a_n\}$ is given as indicated in Proposition \ref{pr:shapeBettidiv}, that is, $a_i = f_i (\prod_{j \neq i} d_j)/d_i$ for all $i \in \{1,\ldots,n\}$, being 
\begin{enumerate}
\item[(a)]  $d_1,\ldots,d_n$ are pairwise prime positive integers, 
\item[(b)] $1 = f_1 = f_2$ and $f_i$ divides $f_{i+1}$ for all $i \in \{2,\ldots,n-1\}$, 
\item[(c)] $\gcd(f_i,d_i) = 1$ for all $i \in \{1,\ldots,n\}$.
\end{enumerate}
Moreover, if we set $D := \prod_{j = 1}^n d_j$ we have that $d_i a_i = f_i D$ for $1 \leq i \leq n$, that $\betti{A} = \{d_i a_i \, \vert \, 2 \leq i \leq n\}$, and  that $d_i a_i$ divides $d_j a_j$  whenever $1 \leq i \leq  j \leq n.$

 The inclusions $\operatorname{Cr}_A  \subseteq \mathcal M_{A} \subseteq \operatorname{Gr}_{A}$ always hold by Proposition \ref{pr:inclusions}, so let us see that $\operatorname{Gr}_A \subseteq \operatorname{Cr}_A$. Take $\x^{\alpha} - \x^{\beta} \in \operatorname{Gr}_{A}$, where $\alpha = (\alpha_1,\ldots,\alpha_n), \beta = (\beta_1,\ldots,\beta_n)\in\mathbb{N}^n$.  We observe that for every $s \in \{1,\ldots,n\}$ we have that  $$(\alpha_s - \beta_s) a_s = \sum_{j \neq s} (\beta_j - \alpha_j) a_j \in \sum_{j \neq s} \Z a_j = \Z \big(\gcd(A -\{a_s\})\big) = \Z d_s$$ Also, for every $s\in\{1,\ldots,n\}$ either $\alpha_s = 0$ or $\beta_s = 0$ (since $\x^{\alpha} - \x^{\beta}$ is primitive), then we get that both $\alpha_s, \beta_s$ are multiples of $ d_s$.
Without loss of generality, we assume now that \begin{equation} \label{eq:rightbig}\ell := \min\{i \, \vert \, \beta_i > 0\} > \min\{i \, \vert \, \alpha_i > 0\}. \end{equation}

{\it Claim 1: }  $\sum_{i = 1}^{t-1} \alpha_i a_i $ is a multiple of $d_t a_t$ for every $t \in \{2,\ldots,\ell\}$.

{\it Proof of claim 1: }  For every $i \geq t$ we have that $d_i a_i$ is a multiple of $d_t a_t$ and, hence, that $\beta_i a_i \equiv 0 \equiv \alpha_i a_i \pmod{d_t a_t}$. As a consequence we have that $\sum_{i = 1}^{t-1} (\alpha_i - \beta_i) a_i \equiv 0 \pmod{d_t a_t}$. By (\ref{eq:rightbig}), we have that $\beta_i = 0$ for all $i \leq t-1 < \ell$ and the claim follows.

{\it Claim 2: } There exist $\alpha_1',\ldots,\alpha_{\ell-1}' \in \mathbb N$ such that $\sum_{i = 1}^{\ell-1} \alpha_i' a_i = d_\ell a_\ell$ and $\alpha_i' \leq \alpha_i$. 

{\it Proof of claim 2: } By (\ref{eq:rightbig}) we have that $\sum_{i = 1}^{\ell-1} \alpha_i a_i \neq 0$ and, by Claim 1, it is a multiple of $d_\ell a_\ell$. Take $s \in \{1,\ldots,\ell-1\}$ such that $\sum_{i = 1}^{s-1} \alpha_i a_i < d_\ell a_\ell \leq  \sum_{i = 1}^{s} \alpha_i a_i$ and define 
\begin{itemize} \item $\alpha_i' = \alpha_i$ for $i = 1,\ldots,s-1$, \item $\alpha_s' :=  \big(d_\ell a_\ell - \sum_{i = 1}^{s-1} \alpha_i a_i)\big)/a_s$, \item $\alpha_{i}' = 0$ for $i = s+1,\ldots,\ell-1$. \end{itemize}
We observe that $\alpha_s' \in \mathbb N$. Indeed, by Claim 1, $\sum_{i = 1}^{s-1} \alpha_i a_i$ is a multiple of $d_s a_s$ and, since $s < \ell$, then $d_\ell a_\ell$ is also a multiple of $d_s a_s$ and, hence, $\alpha_s' = \big(d_\ell a_\ell - \sum_{i = 1}^{s-1} \alpha_i a_i)\big)/a_s \in \mathbb N$. Moreover, by the choice of the $\alpha_i'$ we have that  \[\sum_{i = 1}^{\ell-1} \alpha_i' a_i = \sum_{i = 1}^{s-1} \alpha_i a_i + \alpha_s' a_s = d_\ell a_\ell \] and, by the choice of $s$ we have that $\alpha_s' \leq \alpha_s$.

Now we have that the  binomial $f := \prod_{i = 1}^{\ell-1} x_i^{\alpha_i'} - x_\ell^{d_\ell}$ belongs to $I_A$ and satisfies that  $\prod_{i = 1}^{\ell-1} x_i^{\alpha_i'} \mid \x^{\alpha}$ and $x_\ell^{d_\ell} \mid \x^{\beta}$. But, since $\x^{\alpha} - \x^{\beta} \in \operatorname{Gr}_{A}$, this implies that $f = \x^{\alpha} - \x^{\beta}$. Finally, it suffices to observe that $f \in \operatorname{Cr}_A$ because it is a critical binomial with respect to $x_\ell$. Thus, we have proved that $\operatorname{Cr}_A = \mathcal M_A = \operatorname{Gr}_A$.
 
Since $\mathcal U_A \subseteq \operatorname{Gr}_A = \operatorname{Cr}_A$ we just have to check that all the non-circuits among the critical binomials are not in the universal Gr\"obner basis. Let $g := x_i^{d_i} - \bx^{\alpha} \in \operatorname{Cr}_A$ be a non-circuit. Then the support of $\alpha$ has at least two elements. Let $j,k: 1 \leq j < k \leq n$ be the two smallest values such that both $\alpha_j, \alpha_k$ are non zero. Then, proceeding as before one gets that (1) $d_k$ divides $\alpha_k$, and (2) $\alpha_j a_j \equiv 0 \pmod{d_k a_k}$ which implies that  $d_j f_k/f_j$ divides $\alpha_j$. Now consider $h = x_j^{d_j f_k/f_j} - x_k^{d_k} \in I_A$ and observe that the two monomials in $h$ divide $\bx^{\alpha}$. Hence, for any monomial order $\prec$ we have that $\bx^{\alpha} \in \operatorname{in}_\prec(I_A)$ but it is not a minimal generator of $\operatorname{in}_\prec(I_A)$. It follows that $g$ does not belong to the reduced Gr\"obner basis with respect to $\prec$, which proves that $g \notin \mathcal U_A$.
\end{proof}

Notice that if $\mS$ is a Betti divisible numerical semigroup minimally generated by $A$, then $\mathcal C_A = {\mathcal U}_A = \operatorname{Cr}_A = \mathcal M_A = \operatorname{Gr}_A$ if and only if $I_A$ is generalized robust; equivalently, by \cite[Theorem 4.12]{GT}, $\mS$ has unique Betti degree. Therefore, the inclusion in Theorem \ref{th:bettidivbases} is strict in general.

\subsection{Circuit numerical semigroups}\mbox{}\par

A numerical semigroup $\mS$ with minimal generating set $A = \{a_1,\ldots,a_n\}$ is called a \emph{circuit numerical semigroup} if the toric ideal of $A$ is generated by its set of circuits, that is, if $I_A = \langle \mathcal C_A \rangle$. In this case, we also say that $I_A$ is a \emph{circuit ideal}. The problem of characterizing when a toric ideal is a circuit ideal has been addressed in \cite{BJT,MV} and is widely open.
According to Proposition~\ref{prop:free}, every universally free numerical semigroup is  a circuit numerical semigroup. The converse statement is not true. The following result provides a family of circuit numerical semigroups that contains Betti divisible numerical semigroups and includes some which are not universally free. 

\begin{prop1}\label{pr:circuit}Let $\mS$ be a numerical semigroup minimally generated by $A = \{a_1,\ldots,a_n\}$, where 
$a_j = f_j \prod_{i \neq j} d_i$ with $d_1,\ldots,d_n \in \mathbb Z^+$ pairwise prime and $f_1 = 1, f_2,\ldots,f_n \in \Z^+$ satisfying that $\gcd(d_j,f_j) = 1$ for all $j \in \{1,\ldots,n\}$. Then, $\mS$ is a circuit numerical semigroup and is free for the arrangement $(a_2,\ldots,a_n,a_1)$.
\end{prop1}
\begin{proof} For all $j \in \{2,\ldots,n\}$ we have that 
$\lcm\big(a_j, \gcd(a_1,a_{j+1},\ldots,a_n)\big) =  d_j a_j = f_j d_1 a_1 = \lcm(a_j,a_1).$ Thus, $\mS$ is free for the arrangement $(a_2,\ldots,a_n,a_1).$ Moreover, applying Theorem~\ref{th:free}, we have that
\[ I_A = \langle x_j^{d_j} -  x_1^{f_j d_1} \mid 2 \leq j \leq n \rangle, \]
and we conclude that $\mS$ is a circuit numerical semigroup. 
\end{proof}

As we will see later, for embedding dimension three, all circuit numerical semigroups are of the form of this Proposition. However, as Example \ref{ex:ufnobd} shows, there are universally free numerical semigroups which are not of this form. Furthermore, Example \ref{ex:circ} exhibits a circuit numerical semigroup which is neither of this form nor  universally free.

\begin{ex1}\label{ex:ufnobd}[Continuation of Example~\ref{ex:FAnotBettidiv0}] The numerical semigroup given in Example~\ref{ex:FAnotBettidiv0} is universally free and, hence, a circuit numerical semigroup (and a complete intersection) but it has not the shape described in the statement of Proposition~\ref{pr:circuit}. Indeed, if we write $a_j = f_j \prod_{i \neq j}d_i$ for $i, j \in \{1,2,3,4\}$ with $d_1,\ldots,d_4$ pairwise prime, then neither $11$ nor $13$  divide $d_i$ for all $1 \leq i \leq 4$, $13$ divides $f_1$ and $f_2$, and $11$ divides $f_3$ and $f_4$.
\end{ex1}

\begin{ex1}\label{ex:circ} Consider $\mS = \langle a_1 = 60, a_2 = 280, a_3 = 315, a_4 = 378 \rangle$. We remark that $\mS$ is not free for the arrangement $(a_1,a_2,a_3,a_4)$ because
$\lcm\big(a_1,\gcd(a_2,a_3,a_4)\big) = 7a_1 = 420 \notin \langle a_2,a_3,a_4\rangle$
and, hence, it is not universally free. Moreover, it is free for the arrangement $(a_4,a_3,a_2,a_1)$. Indeed,
\begin{itemize}
    \item $\lcm\big(a_4,\gcd(a_1,a_2,a_3)\big)= 5 a_4 = 6 a_3 \in \langle a_1,a_2,a_3 \rangle,$
     \item $\lcm\big(a_3,\gcd(a_1,a_2)\big)= 4 a_3 = 21 a_1 \in \langle a_1,a_2 \rangle,$ and
     \item $\lcm(a_2,a_1)= 3 a_2 = 14 a_1$.
\end{itemize} Hence, by Proposition \ref{pr:prefree}, we have that: \[I_A = \langle x_1^{14} - x_2^{3},\, x_3^4 - x_1^{21},\, x_4^5 - x_3^6 \rangle. \] 
Moreover, one cannot write the minimal generators of $\mS$ as in Proposition \ref{pr:circuit}.  Indeed, the unique writing of Lemma \ref{lm:uniquewriting} is $d_1 = 7,\, d_2 = 3,\, d_3 = 2,\, d_4 = 5$, $f_1 = f_2 = 4, f_3 = 3$ and $f_4 = 9$.

\begin{verbatim}
gap> s:=NumericalSemigroup(2*3*5*2, 3*5*7*3, 2*3*7*9, 2*5*7*4);
<Numerical semigroup with 4 generators>
gap> IsUniversallyFree(s);
false
gap> IsFree(s);
true
gap> MinimalPresentation(s);
[ [ [ 0, 0, 0, 5 ], [ 0, 0, 6, 0 ] ], [ [ 0, 0, 4, 0 ], [ 7, 3, 0, 0 ] ],
[ [ 0, 3, 0, 0 ], [ 14, 0, 0, 0 ] ] ]
gap> MinimalGenerators(s);
[ 60, 280, 315, 378 ]
gap> AsGluingOfNumericalSemigroups(s);
[ [ [ 60, 280 ], [ 315, 378 ] ], [ [ 60, 280, 315 ], [ 378 ] ],
[ [ 60, 315, 378 ], [ 280 ] ] ]

\end{verbatim}
\end{ex1}

The last result of this subsection characterizes circuit numerical semigroups with embedding dimension three. 

\begin{prop1}\label{pr:em3circuit}
Let $\mS$ be a numerical semigroup with embedding dimension three. Then, $I_A$ is generated by circuits if and only if there exist some pairwise coprime positive integers $d_1,d_2,d_3$, and some $f_2,f_3 \in \Z^+$ such that $\gcd(d_2,f_2) = \gcd(d_3,f_3) =~1$, and $\mS$ is (minimally) generated by $\{d_2 d_3, f_2 d_1 d_3, f_3 d_1 d_2\}$.
\end{prop1}

\begin{proof}
Let $A = \{a_1, a_2, a_3\}$ be the minimal system of generators of $\mS$ and suppose that $I_A$ is generated by circuits. If $\mS$ is not a complete intersection, $I_A$ has a unique Markov basis $M_A$ and the binomials appearing in $M_A$ have full support (see \cite[Section 3]{H}). So, $\mS$ is a complete intersection and, thus, $I_A$ is generated by $M_A = \{x_i^{c_i} - x_j^{c_j}, x_l^{\delta_{kl}} - x_k^{c_k}\}$ where $\delta_{kl} = a_k/\gcd(a_k,a_l), l \in \{i,j\}$ and $\{i,j,k\} = \{1,2,3\}$. For simplicity, let us suppose $i=1, j=2, k=3$ and $l=1$. 

The exponent vectors of the binomials in $M_A$ generate the group $\mathcal L = \big\{(\alpha_1, \alpha_2, \alpha_3) \in \mathbb{Z}^3 \mid \sum_{i=1}^3 \alpha_i a_i = 0\big\}$ and $I_A = \langle \mathbf{x}^{\alpha^+} - \mathbf{x}^{\alpha^-} \mid \alpha \in \mathcal{L}\rangle$ (see \cite[Corollary 3.4]{ST}). Moreover, by \cite[Corollary 2.6]{ES}, $\mathbb{Z}^3/\mathcal{L}$ is torsion free. So, the maximal minors of the matrix \[\left(\begin{array}{ccc} c_1 & -c_2 & 0 \\ \delta_{31} & 0 & -c_3 \end{array}\right),\] name $c_2 c_3, -c_1 c_3$ and $\delta_{31} c_2 $, are $a_1, -a_2$ and $a_3$, respectively; in particular, $c_1$ and $c_2$, $c_2$ and $c_3$, and $\delta_{31}$ and $c_3$ are relatively prime. Now, if $d_1 := \gcd(c_1, \delta_{31})$, then $c_1 = f_2 d_1$ and $\delta_{31} = f_3 d_1$ for some positive intergers $d_2$ and $d_3$ with $\gcd(d_2, d_3) = 1$. Now, taking $d_2 = c_2$ and $d_3 = c_3$ we are done.

The converse is a particular case of Proposition \ref{pr:circuit}.
\end{proof}

Observe that as an immediate consequence of Proposition \ref{pr:em3circuit}, for numerical semigroups with embedding dimension three, we have that: every universally free numerical semigroup is a circuit semigroup, and every circuit semigroup is a complete intersection.

\begin{ex1}\label{ex:circuit-nonFA} In Proposition~\ref{pr:em3circuit}, if one takes $f_2$ and $f_3$ such that there is no divisibility between them, then one has an example of a circuit numerical semigroup that is not Betti divisible (see Proposition~\ref{pr:shapeBettidiv}). For example, taking $d_1 = 2, d_2 = 3,\, d_3 = 5$, $f_2 = 2$ and $f_3 = 3$ one gets $\mS = \langle A \rangle = \langle a_1, a_2, a_3 \rangle$ with $a_1 = d_2 d_3 = 15,\, a_2 = f_2 d_1 d_3 = 20, \, a_3 = f_3 d_1 d_2 = 18$.  One has that $I_A = \langle x_1^{4} - x_2^3,\,x_1^6 - x_3^5 \rangle$
and, hence, it is a circuit ideal and the Betti degrees are $4a_1 = 60$ and $6a_1 = 90$ and, thus, it is not Betti divisible (and, as we will see in the following section, this implies that it is not a universally free numerical semigroup).
\end{ex1}

\begin{quest1} For $A \subseteq \N$, if $I_A$ is a circuit ideal, is $\langle A \rangle$ free? Is $\langle A \rangle$ a complete intersection?
\end{quest1}

\begin{quest1}
For $A\subseteq \N$, is $I_A$ a circuit ideal if  $\betti{A}=\{ \lcm(a,b)\mid a\neq b, a,b\in A\}$?
\end{quest1}

\section{Universally free numerical semigroups with embedding dimension three} \label{sec:ed3}

In this section, we characterize in several ways universally free numerical semigroups with embedding dimension three. An application of the above, is to completely characterize the relation between the size of the toric bases, answering by this way an open question posed in \cite{TTbases}, see Corollary \ref{relation-toric}.

Three-generated numerical semigroups and their  toric ideals have been extensively studied in the literature. Here, we will recall some results concerning them that we will use later; one can find restatements of these results and their proofs in \cite[Chapter~9]{AG} and \cite{H}. 
Let $\mS = \langle a_1,a_2,a_3\rangle$ be a numerical semigroup, then $2 \leq \mu(I_A) \leq 3$ and the Betti degrees of $I_A$ are $\{c_1 a_1, c_2 a_2, c_3a_3\}$. Moreover, $\mu(I_A)  = 2$ or, equivalently, $I_A$ is a complete intersection if and only if there exist $1 \leq i < j \leq 3$ such  that $c_i a_i = c_j a_j$.  Clearly $\mS$ has a unique Betti degree if and only if  $c_1 a_1 = c_2 a_2 = c_3 a_3$.


By Proposition~\ref{pr:inclusions}, we know that $\mathcal{C}_A \cup \operatorname{Cr}_A \subseteq \operatorname{Gr}_A$. The other inclusion does not hold in general as the following example shows.

\begin{ex1}
Let $\mS=\langle 4,6,9\rangle$, which is a complete intersection.
\begin{verbatim}
gap> PrimitiveRelationsOfKernelCongruence([[4],[6],[9]]);
[ [ [ 0, 0, 2 ], [ 0, 3, 0 ] ], [ [ 0, 0, 2 ], [ 3, 1, 0 ] ],
  [ [ 0, 0, 4 ], [ 9, 0, 0 ] ], [ [ 0, 1, 2 ], [ 6, 0, 0 ] ], 
  [ [ 0, 2, 0 ], [ 3, 0, 0 ] ] ]
gap> CircuitsOfKernelCongruence([[4],[6],[9]]);
[ [ [ 3, 0, 0 ], [ 0, 2, 0 ] ], [ [ 9, 0, 0 ], [ 0, 0, 4 ] ], 
  [ [ 0, 3, 0 ], [ 0, 0, 2 ] ] ]
gap> s:=NumericalSemigroup(4,6,9);;
gap> AllMinimalRelationsOfNumericalSemigroup(s);
[ [ [ 0, 3, 0 ], [ 0, 0, 2 ] ], [ [ 3, 0, 0 ], [ 0, 2, 0 ] ], 
  [ [ 3, 1, 0 ], [ 0, 0, 2 ] ] ]    
\end{verbatim}
Notice that $x_1^6-x_2x_3^2$ is a primitive element, but it is neither a circuit nor a critical binomial.

The same holds for the semigroup $\mS=\langle 3,4,5\rangle$, which is not a complete intersection.
\begin{verbatim}
gap> PrimitiveRelationsOfKernelCongruence([[3],[4],[5]]);
[ [ [ 0, 0, 2 ], [ 2, 1, 0 ] ], [ [ 0, 0, 3 ], [ 1, 3, 0 ] ],
  [ [ 0, 0, 3 ], [ 5, 0, 0 ] ], [ [ 0, 0, 4 ], [ 0, 5, 0 ] ],
  [ [ 0, 1, 1 ], [ 3, 0, 0 ] ], [ [ 0, 2, 0 ], [ 1, 0, 1 ] ],
  [ [ 0, 3, 0 ], [ 4, 0, 0 ] ] ]
gap> CircuitsOfKernelCongruence([[3],[4],[5]]);
[ [ [ 4, 0, 0 ], [ 0, 3, 0 ] ], [ [ 5, 0, 0 ], [ 0, 0, 3 ] ],
  [ [ 0, 5, 0 ], [ 0, 0, 4 ] ] ]
gap> s:=NumericalSemigroup(3,4,5);
<Numerical semigroup with 3 generators>
gap> AllMinimalRelationsOfNumericalSemigroup(s);
[ [ [ 1, 0, 1 ], [ 0, 2, 0 ] ], [ [ 2, 1, 0 ], [ 0, 0, 2 ] ],
  [ [ 3, 0, 0 ], [ 0, 1, 1 ] ] ]
\end{verbatim}
In this case, $x_3^3-x_1x_2^3$ is primitive, but it is neither a circuit nor a critical binomial.

\end{ex1}

Next, we prove the main theorem of this section, in which we characterize the universally free numerical semigroups with embedding dimension three.

\begin{thm1}\label{th:char-3-FA}
Let $\mS$ be a numerical semigroup minmally generated by $A = \{a_1, a_2, a_3\}$. The following are equivalent:
\begin{itemize}
\item[(a)] $\mS$ is a universally free.
\item[(b)] $\mS$ is Betti divisible.
\item[(c)] There exist positive integers $d_1$, $d_2$, $d_3$ and $f_3$ such that $A = \{d_2d_3, d_1d_3, f_3 d_1d_2\}$ and $\betti{A}=\{d_1d_2d_3,f_3 d_1d_2d_3\}$; in particular, $\gcd(d_i,d_j) = 1$ for $ i \neq j$, and $\gcd(f_3,d_3) =~1$.
\end{itemize}
Moreover, up to permutation of the indeterminates, in this case we have that
\begin{align*}
\mathcal{M}_A & =\big\{x_1^{d_1}-x_2^{d_2}, x_3^{d_3}-x_1^{f_3 d_1},x_3^{d_3}-x_1^{(f_3-1)d_1}x_2^{d_2},\dots,x_3^{d_3}-x_1^{d_1}x_2^{(f_3-1)d_2},x_3^{d_3}-x_2^{f_3 d_2}\big\} \\ & =\operatorname{Cr}_A 
=\operatorname{Gr}_A,
\end{align*}
and $\mathcal U_A = \big\{x_1^{d_1} - x_2^{d_2},\, x_1^{f_3 d_1} - x_3^{d_3},\, x_2^{f_3 d_2} - x_3^{d_3}\big\} = \mathcal C_A$.
\end{thm1}

\begin{proof}
Write $\delta_{ij}=a_i/\gcd(a_i,a_j)$, $i\neq j$, and let $c_i$ be as in \eqref{eq:cb}. Recall that $\betti{A}=\{c_1a_1,c_2a_2,c_3a_3\}$.

First, assume that $\mS$ is universally free. As $\mS$ is a complete intersection, the cardinality of $\betti{A}$ is less than three. If $\mS$ has a unique Betti degree, then the result follows directly from \cite[Theorem~12]{GOR}. Thus, suppose that the cardinality of $\betti{A}$ is two. Without loss of generality, let us assume $c_1a_1=c_2a_2\neq c_3a_3$. So, by \cite[Proposition 2.3]{KO} and Proposition \ref{prop:free2}, we have that both $\{x_1^{c_1} - x_2^{c_2}, x_3^{c_3} - x_1^{\delta_{31}}\}$ and $\{x_1^{c_1} - x_2^{c_2}, x_3^{c_3} - x_2^{\delta_{32}}\}$ are a Markov bases of $I_A$; in particular, 
$x_1^{\delta_{31}} - x_2^{\delta_{32}} \in I_A$. Therefore, $c_1$ divides $\delta_{31}$ and $c_2$ divides $\delta_{32}$. Hence, $c_1 a_1 = c_2 a_2$ divides $\delta_{31} a_1 = \delta_{32} a_2 = a_3 c_3$, that is, $\mS$ is Betti divisible.

The equivalence of conditions (b) and (c) is nothing more than Proposition \ref{pr:shapeBettidiv} in the embedding dimension three case. 

Finally, if condition (c) holds, then one can easily check that $\mS$ is free for all possible arrangements of $\{a_1,a_2,a_3\}$ and the first part of the proof is completed.

Now, if $\mS$ has a unique Betti degree, then by \cite[Theorem 6]{GOR} the following toric bases are equal: $\operatorname{Cr}_A = \mathcal M_A = \operatorname{Gr}_A = \mathcal C_A = {\mathcal U}_A.$
If $\mS$ has two Betti degrees (or, equivalently, $f_3 > 1$), Theorem \ref{th:bettidivbases} guarantees that $\operatorname{Cr}_A = \mathcal M_A = \operatorname{Gr}_A$ and $\mathcal C_A = {\mathcal U}_A.$
So, it only remains to show that  
\[
\mathcal M_A=\big\{x_1^{d_1}-x_2^{d_2}, x_3^{d_3}-x_1^{f_3 d_1},x_3^{d_3}-x_1^{(f_3-1)d_1}x_2^{d_2},\dots,x_3^{d_3}-x_2^{f_3 d_2}\big\}.
\]
Since, by \cite[Proposition 2.3]{KO}, $x_1^{c_1}, x_2^{c_2}$ and $x_3^{c_3}$ have to appear in the binomials in any Markov basis of $I_A$, therefore the Markov bases of $I_A$ are of the form $\{x_1^{c_1} - x_2^{c_2}, x_3^{c_3} - x_1^{\alpha_1} x_2^{\alpha_2}\}$ for some $(\alpha_1, \alpha_2) \in \mathbb{N}^2$. To determine the possible $(\alpha_1, \alpha_2)$ one has to obtain all the expressions of $c_3a_3$ in terms of $a_1$ and $a_2$. Take $A'=\{a_1,a_2\}$. Note that $c_i = d_i,\ i = 1, 2, 3$, and that $c_3a_3=f_3 d_1d_2d_3=f_3 d_1a_1$, and so it suffices to see all the different expressions of $f_3 d_1a_1$ in the monoid generated by $A'$. The only generator of $I_{A'}$ is precisely $x_1^{d_1}-x_2^{d_2}$, and so the expressions of $f_3 d_1a_1$ are $\{f_3 d_1a_1,(f_3-1)d_1a_1+d_2a_2,\dots, d_1a_1+(f_3-1)d_2a_2,f_3 d_2a_2\}$ and we are done.
\end{proof}

\begin{ex1}
Let $\mS=\langle 10, 15, 18\rangle$. In this case, $d_1=3$, $d_2=2$, $d_3=5$ and $f_3=3$.
\begin{verbatim}
gap> s:=NumericalSemigroup(10,15,18);;
gap> IsUniversallyFree(s);
true
gap> BettiElements(s);
[ 30, 90 ]
gap> AllMinimalRelationsOfNumericalSemigroup(s);
[ [ [ 0, 6, 0 ], [ 0, 0, 5 ] ], [ [ 3, 0, 0 ], [ 0, 2, 0 ] ],
[ [ 3, 4, 0 ], [ 0, 0, 5 ] ], [ [ 6, 2, 0 ], [ 0, 0, 5 ] ], 
[ [ 9, 0, 0 ], [ 0, 0, 5 ] ] ]
gap> PrimitiveRelationsOfKernelCongruence([[10],[15],[18]]);
[ [ [ 0, 0, 5 ], [ 0, 6, 0 ] ], [ [ 0, 0, 5 ], [ 3, 4, 0 ] ],
[ [ 0, 0, 5 ], [ 6, 2, 0 ] ], [ [ 0, 0, 5 ], [ 9, 0, 0 ] ],
[ [ 0, 2, 0 ], [ 3, 0, 0 ] ] ]
\end{verbatim}

And, as the following extract from SageMath \cite{sage} shows, the universal Gr\"obner basis of $I_A$ equals the set of circuits.
\begin{verbatim}
R.<x,y,z> = QQ[]
I = R.ideal([x^3 - y^2, y^6 - z^5])
G = I.groebner_fan()
G.reduced_groebner_bases()
[[x^3 - y^2, y^6 - z^5],
 [x^9 - z^5, -x^3 + y^2],
 [-x^3 + y^2, -x^9 + z^5],
 [x^3 - y^2, -y^6 + z^5]]
\end{verbatim}

\end{ex1}

In \cite{TTbases} the authors study the relative size of toric bases. In particular they prove in \cite[Theorem 3.1]{TTbases} that the size of the elements of the Graver basis, the universal Gr\"obner basis and the set of the circuits of a toric ideal cannot be bounded above by a polynomial on the size of a Markov basis of $I_A$. In Section 4 of the same paper they leave as an open problem to prove that the size of the Graver basis cannot be bounded above by a polynomial expression on the size of the universal Gr\"obner basis or the set of circuits (see \cite[Figure 2]{TTbases} for more details). Here, in Theorem \ref{th:char-3-FA}, we provide a family of examples where $\mathcal{U}_A = \mathcal{C}_A$ has size three and $\operatorname{Gr}_A = \mathcal{M}_A$ is arbitrarily large. In particular, it follows the next corollary.

\begin{cor1}\label{relation-toric} The size of the Graver basis $\operatorname{Gr}_A$ or the universal Markov basis $\mathcal M_A$ of a toric ideal $I_A$, cannot be bounded above by a function on the size of the universal Gr\"obner basis $\mathcal U_A$ or the set of the circuits $\mathcal C_A$ of $I_A$.
\end{cor1}

\section{Conclusion / Arbitrary embedding dimension}

Consider the following families of numerical semigroups:
\begin{itemize}
\item[$\mathcal F_0$:] Betti divisible numerical semigroups.
\item[$\mathcal F_1$:] Numerical semigroups where $\mathcal M_A = \operatorname{Gr}_A$.
\item[$\mathcal F_2$:] Numerical semigroups where $\mathcal C_A = \mathcal U_A$.
\item[$\mathcal F_3$:] Universally free semigroups.
\item[$\mathcal F_4$:] Numerical semigroups where $\mathcal C_A \subseteq \mathcal M_A$.
\item[$\mathcal F_5$:] Circuit numerical semigroups.
\end{itemize}

In the previous sections we proved that whenever $\mS = \langle A \rangle$ with $A = \{a_1,a_2,a_3\}$ is a numerical semigroup of embedding dimension $3$, we have that $\mathcal F_0 = \mathcal F_1 = \mathcal F_2 = \mathcal F_3 = \mathcal F_4 \subsetneq \mathcal F_5$. In this section we study the inclusion  relations that hold between these families of numerical semigroups in arbitrary embedding dimension.


A first observation is that the families $\mathcal F_i$ for $i \in \{0,1,2,3,4\}$ are \emph{closed under elimination of variables}, that is, if $\mS = \langle A \rangle$ belongs to $\mathcal F_i$ and $A' \subseteq A$, then $\mS' = \langle A' \rangle$ also belongs to $\mathcal F_i$.  Indeed, this is a consequence of Proposition~\ref{pr:shapeBettidiv} (for $i = 0$), the definition of universally free semigroups (for $i = 3$) and Propositions~\ref{pr:elimination} and \ref{genrobustfewvariables} (for $i = 1,2,4$).
However, as the following example shows, the family of circuit numerical semigroups is not closed under elimination of variables. 
\begin{ex1} Consider $A = \{a_1,a_2,a_3,a_4\}$ with $a_1 = 30,\, a_2 = 36,\, a_3 = 40,\, a_4 = 75$, then $I_A = \langle x_2^5 - x_1^6,\, x_3^3 - x_1^4,\, x_4^2 - x_1^5 \rangle$ is a circuit ideal. However, taking $A' = \{a_2,a_3,a_4\}$ we have that $I_{A'}$ is not a circuit ideal because $300$ is a Betti degree but is not the $A$-degree of any circuit. 
\begin{verbatim}
gap> s:=NumericalSemigroup(36,40,75);
<Numerical semigroup with 3 generators>
gap> A:=MinimalGenerators(s);
[ 36, 40, 75 ]
gap> BettiElements(s);
[ 300, 360 ]
gap> Set(Combinations(A,2),Lcm);
[ 360, 600, 900 ]
\end{verbatim}
\end{ex1}

Next, we are going to show that the numerical semigroup of Example \ref{ex:FAnotBettidiv0}, which we know that is not Betti divisible, belongs to  $\mathcal F_1 \cap \mathcal F_2$.

\begin{ex1} \label{ex:nobetidivequalbases} [Continuation of Example \ref{ex:FAnotBettidiv0}]
We know that the numerical semigroup $\mS = \langle a_1 = 390, a_2 = 546, a_ 3 = 770, a_4 = 1155 \rangle$ of Example~\ref{ex:FAnotBettidiv0} is universally free and not Betti divisible. Let us check that $\mS$ belongs to both $\mathcal F_1$ and to $\mathcal F_2$.

The following GAP code checks that $\mathcal M_A = \operatorname{Gr}_A$ (and also that both sets consist of 170 elements and $x_1^7 x_2^{50} - x_3^3 x_4^{24}$ belongs to them).
\begin{verbatim}
gap>  s:=NumericalSemigroup(390,546,770,1155);
<Numerical semigroup with 4 generators>
gap> A:=AllMinimalRelationsOfNumericalSemigroup(s);;
gap> B:=PrimitiveRelationsOfKernelCongruence(List(MinimalGenerators(s),g->[g]));;
gap> Set(A,SortedList)=Set(B,SortedList);
true
gap> Size(A);
170
gap> IsMinimalRelationOfNumericalSemigroup([[7,50,0,0],[0,0,3,24]],s);
true
\end{verbatim}

With SageMath \cite{sage} one can compute the set of all reduced Gr\"obner bases of $I_A$.
\begin{verbatim}
R.<x,y,z,t> = QQ[]
I = R.ideal([t^2-z^3, x^7-y^5, t^26 - y^55])
G = I.groebner_fan()
G.reduced_groebner_bases()
[[z^3 - t^2, x^7 - y^5, y^55 - t^26],
 [x^77 - t^26, z^3 - t^2, -x^7 + y^5],
 [z^3 - t^2, -x^7 + y^5, -x^77 + t^26],
 [y^55 - z^39, x^7 - y^5, -z^3 + t^2],
 [x^77 - z^39, -z^3 + t^2, -x^7 + y^5],
 [-z^3 + t^2, -x^7 + y^5, -x^77 + z^39],
 [-z^3 + t^2, x^7 - y^5, -y^55 + z^39],
 [z^3 - t^2, x^7 - y^5, -y^55 + t^26]]
\end{verbatim}
Hence, one has that $\mathcal U_A = \mathcal C_A$, and $\mS$ belongs to $\mathcal F_2$.

\medskip

The following result proves the inclusion relations depicted in Figure \ref{fig:inclusions}
\begin{prop1}\label{pr:inclusions2}The following inclusions hold:
\begin{itemize}
    \item[(a)] $\mathcal F_0 \subsetneq \mathcal F_1 \cap \mathcal F_2,$
    \item[(b)] $\mathcal F_1 \cup \mathcal F_3 \subseteq \mathcal F_4,$
    \item[(c)] $\mathcal F_2  \subseteq \mathcal F_3 \subsetneq \mathcal F_5.$
\end{itemize}
\end{prop1}

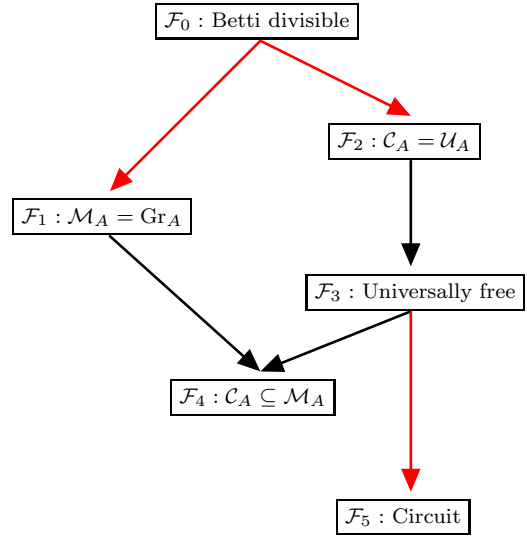
\begin{figure}
\begin{tikzpicture}[line cap=round,line join=round,>=triangle 45,x=1cm,y=1cm]
\clip(0,1) rectangle (9,9);
\draw [->,line width=1pt,color=red] (4,8) -- (6,7);
\draw [->,line width=1pt,color=red] (4,8) -- (2,6);
\draw [->,line width=1pt] (6,6.4) -- (6,5);
\draw [->,line width=1pt] (2,5.4) -- (4,3.6);
\draw [->,line width=1pt] (6,4.4) -- (4,3.6);
\draw [->,line width=1pt,color=red] (6,4.4) -- (6,2);
\begin{scriptsize}
\draw (2.5,8.6) node[anchor=north west] {\framebox{$\mathcal{F}_0:$ Betti divisible}};
\draw (4.8,7) node[anchor=north west] {\framebox{$\mathcal{F}_2: \mathcal{C}_A = \mathcal{U}_A$}};
\draw (0.6,6) node[anchor=north west] {\framebox{$\mathcal{F}_1:\mathcal{M}_A = \operatorname{Gr}_A$}};
\draw (4.5,5) node[anchor=north west] {\framebox{$\mathcal{F}_3:$ Universally free}};
\draw (2.7,3.6) node[anchor=north west] {\framebox{$\mathcal{F}_4: \mathcal{C}_A \subseteq \mathcal{M}_A$}};
\draw (4.9,2) node[anchor=north west] {\framebox{$\mathcal{F}_5:$ Circuit}};
\end{scriptsize}
\end{tikzpicture}
\caption{Arrows show inclusions between classes of numerical semigroups. Red arrows indicate that inclusion is strict.} \label{fig:inclusions}
\end{figure}

\begin{proof}We have that $\mathcal F_0 \subseteq \mathcal F_1 \cap \mathcal F_2$ by Theorem \ref{th:bettidivbases} and the inclusion is strict by Example \ref{ex:nobetidivequalbases}. 

The inclusion $\mathcal F_3 \subseteq \mathcal F_4$ follows from Proposition \ref{prop:free2}, and $\mathcal F_1 \subseteq \mathcal F_4$ because $\mathcal C_A \subseteq \operatorname{Gr}_A = \mathcal M_A$ (see Proposition \ref{pr:inclusions}). Hence, $\mathcal F_1 \cup \mathcal F_3 \subseteq \mathcal F_4$. 

Also, $\mathcal F_3 \subseteq \mathcal F_5$ by Proposition \ref{prop:free}.(c). Inclusion is strict by Example \ref{ex:circ}.

Hence, it only remains to prove that $\mathcal F_2 \subseteq \mathcal F_3$. Consider a numerical semigroup $\mS = \langle A \rangle$ such that $\mathcal C_A = \mathcal U_A$ and $\prec$  a monomial order. If follows that $\mathcal U_A$ is a Gr\"obner basis of $I_A$ with respect to $\prec$. Since $\mathcal C_A = \mathcal U_A$ we derive that $\operatorname{in}_\prec(I_A)$ is generated by monomials $\textbf{x}^{\textbf{u}}$ whose support has only one element (i.e., by powers of variables). As $\operatorname{ht}(I_A) = n-1$ it follows that 
there exists a $j \in \{1,\ldots,n\}$ and some $e_1,\ldots,e_{j-1},e_{j+1},\ldots,e_n \in \mathbb Z^+$ such that $\operatorname{in}_\prec(I_A) = \langle x_1^{e_1},\ldots,x_{j-1}^{e_{j-1}},x_{j+1}^{e_{j+1}},\ldots,x_n^{e_n} \rangle$ and, thus, the reduced Gr\"obner basis of  $I_A$ with respect to $\prec$ has $n-1$ elements. By  Theorem \ref{universally-Grobner}, this implies that $\mS$ is universally free.  
\end{proof}

Proposition~\ref{pr:inclusions} and Proposition~\ref{prop:free2} imply that $\mathcal{C}_A \cup \operatorname{Cr}_A \subseteq \mathcal M_A$ for universally free numerical semigroups. For $n > 3$, the opposite inclusion does not always hold. Again the numerical semigroup $\mS = \langle 390, 546, 770,  1155 \rangle$ serves as a counterexample because the binomial 
$x_1^7 x_2^{50} - x_3^3 x_4^{24}$ belongs to $\mathcal M_A$, but is not a circuit or a critical binomial. 

A toric ideal $I_A$ is said to be {\it critical} if it is generated by the set of critical binomials, i.e., if $I_A = \langle \operatorname{Cr}_A \rangle$ (see, e.g., \cite{KO}) . As a consequence of Theorem \ref{th:bettidivbases}, every Betti divisible numerical semigroup defines a critical toric ideal (because $\operatorname{Cr}_A = \mathcal M_A$). Thus is no longer true for universally free numerical semigroups and, again, $\mS = \langle 390,546,770,1155\rangle$ serves as a counterexample. In this example there are only two critical binomials $\operatorname{Cr}_A = \{x_1^7 - x_2^5,\, x_3^3 - x_4^2\}$ and, hence, it is impossible that they generate the height $3$ of the ideal $I_A$.
\end{ex1}

 
 

In Proposition \ref{pr:inclusions2} we proved several inclusions among the classes $\mathcal F_i$, moreover we wonder if the following statement (which includes Conjecture \ref{cj:conj1}), holds.

\begin{conj1} The families $\mathcal F_i$ equal for $i = 1,2,3,4$.
\end{conj1}

In Proposition \ref{pr:inclusions2} we proved that $\mathcal F_1 \subseteq \mathcal F_4$ and $\mathcal F_2 \subseteq \mathcal F_4$, and computational experiments support that these three classes might be equal. These results seems to suggest that, at least in the context of toric ideals associated to numerical semigroups, there are deeper connections between the toric bases apart from those described in Proposition \ref{pr:inclusions}. We believe that these relations among toric bases of general toric ideals deserve further studies.

Finally, we believe that it would be interesting to characterize universally free numerical semigroups in terms of gluings. We propose the following conjecture. 

\begin{conj1}\label{conj-glue}
Let $\mS$ be a numerical semigroup minimally generated by $A = \{a_1,\ldots,a_n\}$. Then, $\mS$ is universally free if and only if there exists a nontrivial partition of $A$ in $\{B,C\}$ in such a way that $\lcm\big(\gcd(B),\gcd(C)\big) = \lcm(A)$ and both $\langle B \rangle$ and $\langle C \rangle$ are universally free numerical semigroups. 
\end{conj1}

\begin{ex1}\label{last}[Continuation of Example~\ref{ex:FAnotBettidiv0}]
We remind that $\mS = \langle a_1 = 390, a_2 = 546, a_ 3 = 770, a_4 = 1155 \rangle$. Taking $B = \{a_1,a_2\}$ and $C = \{a_3,a_4\}$ we have that \[ \lcm\big(\gcd(B), \gcd(C)\big) = \lcm(2 \cdot 3 \cdot 13,\ 5 \cdot 7 \cdot 11) =  2 \cdot 3 \cdot 5 \cdot 7 \cdot 11 \cdot 13 = \lcm(A).\] 

\begin{verbatim}
gap> AsGluingOfNumericalSemigroups(s);
[ [ [ 390 ], [ 546, 770, 1155 ] ], [ [ 390, 546 ], [ 770, 1155 ] ],
  [ [ 390, 546, 770 ], [ 1155 ] ], [ [ 390, 546, 1155 ], [ 770 ] ],
  [ [ 390, 770, 1155 ], [ 546 ] ] ]
\end{verbatim}
\end{ex1}

By using the \texttt{GAP} package \texttt{numericalsgps} we have been able to check that Conjecture~\ref{conj-glue} holds for semigroups generated by four elements ranging from 10 to 500. We believe that a proof of this conjecture would shed light on the structure of universally free numerical semigroups. To the best of our knowledge, there is no direct link between this conjecture and the previously stated ones. 

\section*{Acknowledgments}

The first and the fourth authors are partially supported by the Spanish MICINN ALCOIN (PID2019-104844GB-I00) and by the ULL funded research project MACACO.

The second and third authors are partially supported by the Proyecto de Excelencia de la Junta de Andalucía (ProyExcel\_00868). The second author acknowledges financial support from the Spanish Ministry of Science and Innovation (MICINN), through the “Severo Ochoa and María de Maeztu Programme for Centres and Unities of Excellence” (CEX2020-001105-M). The third author is partially supported by the Proyecto de investigación del Plan Propio - UCA 2022-2023 (PR2022-011).


\end{document}